\documentclass[twoside,a4paper,reqno]{amsart}

\usepackage{style}
\usepackage{graphicx}

\usepackage{pdfsync}

\usepackage{xcolor}
\definecolor{webgreen}{rgb}{0,.5,0}
\definecolor{webbrown}{rgb}{.6,0,0}
\definecolor{RoyalBlue}{cmyk}{1, 0.50, 0, 0}
\usepackage[colorlinks=true, breaklinks=true, urlcolor=webbrown, linkcolor=RoyalBlue, citecolor=webgreen,backref=page]{hyperref}

\usepackage{mathrsfs}

\usepackage{courier} 
\normalfont
\usepackage{epsfig, graphicx, subfigure}
\usepackage{verbatim, setspace}
\usepackage{amsmath, amssymb}

\numberwithin{equation}{section}

\usepackage{courier} 
\normalfont

\def\cal{\mathcal}
\let\Re=\undefined
\DeclareMathOperator{\Re}{Re}
\let\Im=\undefined
\DeclareMathOperator{\Im}{Im}

\def\ge{\geqslant}
\def\le{\leqslant}

\begin{document}

\title[NLS]{Sobolev norms of $L^2-$ solutions to NLS}
  \author{Roman V.~Bessonov, \quad  Sergey A.~Denisov}

\address{
\begin{flushleft}
Roman Bessonov: bessonov@pdmi.ras.ru\\\vspace{0.1cm}
St.\,Petersburg State University\\  
Universitetskaya nab. 7-9, 199034 St.\,Petersburg, RUSSIA\\
\vspace{0.1cm}
St.\,Petersburg Department of Steklov Mathematical Institute\\ Russian Academy of Sciences\\
Fontanka 27, 191023 St.Petersburg,  RUSSIA\\
\end{flushleft}
}
\address{
\begin{flushleft}
Sergey Denisov: denissov@wisc.edu\\\vspace{0.1cm}
University of Wisconsin--Madison\\  Department of Mathematics\\
480 Lincoln Dr., Madison, WI, 53706, USA\\
\end{flushleft}
}

\thanks{The work of RB in Sections 2 and 3 is supported by the Russian Science Foundation grant 19-71-30002. The work of SD  {in the rest of the paper} is supported by the grant NSF DMS-2054465  and Van Vleck Professorship Research Award. RB is a Young Russian Mathematics award winner and would like to thank its sponsors and jury.}

\begin{abstract} 
We apply inverse spectral theory to study Sobolev norms of solutions to the nonlinear Schr\"{o}dinger equation. For initial datum $q_0\in L^2(\R)$ and $s\in [-1,0]$, 
we prove that there exists a conserved quantity which is equivalent to $H^s(\R)$-norm of the solution. 
\end{abstract} 

\vspace{1cm}


\maketitle
\setcounter{tocdepth}{1}




\section{Introduction}
\noindent The classical defocusing  nonlinear Schr\"{o}dinger equation (NLS)  \cite{FTbook,TaoNDE,ZS71} on the real line has the form
\begin{equation}\label{nls}
\begin{cases}
i \frac{\partial q}{\partial t} = -\frac{\partial^2 q}{\partial \xi^2} + 2|q|^2q, \\
q\big\rvert_{t = 0} = q_0,
\end{cases}
\qquad \xi \in \R, \quad t \in \R.
\end{equation}
It is known that for sufficiently regular initial datum $q_0$ the unique classical solution $q = q(\xi, t)$ exists globally in time. For example, if $q_0$ lies in the Schwartz class $\Sch(\R)$, then $q(\cdot, t)\in \Sch(\R)$ for all $t\in \R$. The long-time asymptotics of $q$ is  known \cite{DZ2003,ZM1976,DZ1994}. For less regular initial datum $q_0$, one can define the solution  by an approximation argument (see, e.g., \cite{yap}):
\begin{Thm}\label{yap}
Let $q_{0} \in L^2(\R)$, and let $q_{0,n} \in \Sch(\R)$ converge to $q_0$ in $L^2(\R)$. Denote by $q_n(\xi,t)$ the solution of \eqref{nls} corresponding to $q_{0,n}$. We have
$$
\lim_{n \to +\infty}\|q_n(\cdot, t)-q(\cdot, t)\|_{L^2(\R)} = 0, \qquad t \in \R,
$$ 
for some function $q(\xi, t): \R^{2} \to \R$ that does not depend on the choice of the sequence $q_{0,n}$.
\end{Thm}  
The function $q$ in Theorem \ref{yap} is called the $L^2$--solution of \eqref{nls} corresponding to the initial datum $q_0 \in L^2(\R)$. It is clear that such a solution is unique.  The total  energy of the solution is its  $L^2(\R)$-norm and it is conserved in time:
$$
\|q(\cdot, t)\|_{L^2(\R)} = \|q_0\|_{L^2(\R)}, \qquad t \in \R\,.
$$   
By Plancherel's formula, it is equal to
$\|(\cal{F}q)(\cdot, t)\|_{L^2(\R)}$
where $\cal F$ stands for the Fourier transform.
In this paper, we work with Sobolev spaces $H^s(\R)$, $s\in \R$. The $H^s(\R)$-norm of a function $f\in \Sch(\R)$ is defined by 
\begin{equation}\label{n1}
\|f\|_{H^s(\R)}=\left(\int_{\R}(1+|\eta|^2)^{s}|(\cal{F}f)(\eta)|^2d\eta\right)^{\frac 12}\,.
\end{equation}
The space $H^s(\R)$ is the completion of $\Sch(\R)$ with respect to this norm. Equivalently, one can define it by
\[
H^s(\R)=\{f\in \Sch'(\R): (1+|\eta|^2)^{{\frac s2}}\cal{F}f\in L^2(\R)\}\,,
\]
where $\Sch'(\R)$ is the space of tempered distribution. 

\medskip

In contrast to the linear Schr\"odinger equation for which all Sobolev norms are conserved, the solutions of NLS can exhibit
inflation of Sobolev norm $H^{s}(\R)$ for $s\le -\frac 12$ (see, e.g., \cite{christ,kishimoto} for details). Specifically, given an arbitrarily small positive $\eps$ and $s\le-\frac 12$, there exists a solution $q$ to \eqref{nls} that satisfies
\begin{equation}\label{tao}
q_0\in \Sch(\R), \qquad \|q_0\|_{H^s(\R)}\le \eps, \qquad \|q(\cdot,\eps)\|_{H^s(\R)}\ge \eps^{-1}\,,
\end{equation}
see \cite{christ} for that construction.
This result is related to the ``high-to-low frequency cascade''. It occurs when for  initial datum $q_0\in \Sch(\R)$, a part of  $L^2(\R)$-norm of $q$, when written on the Fourier side,
moves from high to low frequencies as time increases. The Sobolev norms with negative index $s$ can be used to capture this phenomenon. Indeed, since $\|q(\cdot,t)\|_{L^2(\R)}$ is time-invariant and the weight $(1+\eta^2)^s$ in \eqref{n1} vanishes at infinity when $s<0$,  the transfer of $L^2$-norm from high to low values of frequency $\eta$ makes the $H^s(\R)$-norm grow. 

\smallskip

For NLS, the  inflation of $H^s(\R)$-norm can not happen for $s>-{\frac 12}$.
In \cite{tataru}, Koch and Tataru  discovered the set of conserved quantities which agree with $H^{s}(\R)$-norm up to a quadratic term for a small value of $\|q_0\|_{H^s(\R)}$ and $s>-\frac 12$. As a corollary, they obtained the bounds on $\|q(\cdot,t)\|_{H^s(\R)}$ that are uniform in time:
\begin{equation}\label{koch}
\|q(\cdot,t)\|_{H^s(\R)}\le C(s)
\left\{\begin{matrix}
\mathcal{R}+\mathcal{R}^{1+2s}, & s>0,\\
\mathcal{R}+\mathcal{R}^{\frac{1+4s}{1+2s}}, & s\in (-\frac 12,0)\,,
\end{matrix}\right. \quad \mathcal{R}=\|q_0\|_{H^s(\R)}\,.
\end{equation}
 {In  \cite{Killip18}, Killip, Vi\c{s}an, and Zhang proved a similar estimate using a different method. }The estimates on the growth of $H^s(\R)$-norms are related to questions of well-posedness and ill-posedness of NLS in Sobolev classes which have been extensively studied previously, see, e.g., \cite{kenig,tataru,kishimoto, BC20, CK17, CJT03, OS12}.\smallskip

In our paper, we use some recent results in the inverse spectral theory \cite{BD2017,BD2019,BD2022} to show that there are conserved quantities of NLS which agree with $H^{s}(\R)$-norm provided that $s\in [-1,0]$ and the value of  $\|q_0\|_{L^2(\R)}$ is under control. We apply our analysis to prove the following theorem.
\begin{Thm}\label{t5}
Let $q_0 \in L^2(\R)$  and let $q = q(\xi, t)$ be the solution of \eqref{nls} corresponding to $q_0$. Then,  
\begin{equation}\label{eq19}
C_1(1+\|q_0\|_{L^2(\R)})^{2s} \|q_0\|_{H^{s}(\R)}\le \|q(\cdot, t)\|_{H^{s}(\R)} \le C_2(1+\|q_0\|_{L^2(\R)})^{-2s} \|q_0\|_{H^{s}(\R)},
\end{equation} 
where $\,t \in \R\,, \,s\in [-1,0]$, and $C_1$ and $C_2$ are two positive absolute constants.
\end{Thm}
This result shows, in particular, that for a given function $q_0: \|q_0\|_{L^2(\R)}=1$ whose $L^2(\R)$-norm is concentrated on high frequencies, we will never see a significant part of  $L^2(\R)$-norm of the solution $q$ moving to the low frequencies. That limits the ``high-to-low frequency cascade'' we discussed above.  The close inspection of construction used in \cite{christ} shows that the function $q_0$ in \eqref{tao}  has $H^s(\R)$-norm  smaller than  {$\eps$} but its $L^2(\R)$-norm is  large when  {$\eps$} is small. Hence, the bounds in Theorem~\ref{t5} do not contradict the estimates in \eqref{tao} when $s\in [-1,-\frac 12]$.  We do not know whether Theorem~\ref{t5}  holds for $s<-1$.\smallskip

The main idea of the proof of Theorem~\ref{t5} is based on the analysis of the conserved quantity $a(z)$,  $\Im z > 0$, which is a coefficient in the transition matrix for the Dirac equation with potential $q=q(\cdot,t)$. We take $z=i$ and show that $\log |a(i)|$ is related to a certain quantity $\widetilde\K_{Q}$ (see the Lemma~\ref{l88} below) that characterizes both size and oscillation of $q$. Using $\widetilde\K_{Q}$ in the context of NLS is the main novelty of our work. 
We study $\widetilde\K_{Q}$ and show that it is equivalent to $H^{-1}(\R)$ norm of $q$ with constants that depend on its $L^2(\R)$-norm.  That gives the estimate \eqref{eq19} for $s=-1$ and  the intermediate range of $s\in (-1,0)$ is handled by interpolation. Our analysis relies heavily on the recent results \cite{BD2017,BD2019,BD2022} that characterize Krein -- de Branges canonical systems and the Dirac operators whose spectral measures belong to the Szeg\H{o} class on the real line. We also establish the framework that allows working with NLS in the context of well-studied Krein systems.

\smallskip

 \bigskip

\noindent {\bf Notation}\medskip
\begin{itemize}
\item The symbol $I$ stands for $2\times 2$ identity matrix $I=\idm$ and symbol $J$ stands for $J=\left(\begin{smallmatrix} 0& -1\\1 &0\end{smallmatrix} \right)$.  {Constant matrices $\sigma_{3}$, $\sigma_{\pm}$, $\sigma$ are defined in \eqref{cmatrices}.}
\item  {For} a measurable set $S\subset \R$, we say that $f\in L^1_{\rm \loc}(S)$ is $f\in L^1(K)$ for every compact $K\subset S$.
\item The Fourier transform of a function $f$ is defined by
$$
(\F f)(\eta) = \frac{1}{\sqrt{2\pi}}\int_{\R}f(x)e^{-i\eta x}\,dx.
$$    

\item The symbol $C$, unless we specify explicitly, denotes the absolute constant which can change the value from formula to formula. If we write, e.g., $C(\alpha)$, this defines a positive function of parameter $\alpha$.

\item For two non-negative functions $f_1$ and $f_2$, we write $f_1\lesssim f_2$ if  there is an absolute
constant $C$ such that $f_1\le Cf_2$ for all values of the arguments of $f_1$ and $f_2$. We define $\gtrsim$
similarly and say that $f_1\sim f_2$ if $f_1\lesssim f_2$ and
$f_2\lesssim f_1$ simultaneously. If $|f_3|\lesssim f_4$, we will write $f_3=O(f_4)$.

\item Symbols $\{e_j\}$ are reserved for the standard basis in $\C^2$: $e_1=\oz$, $e_2=\zo$.

\item For matrix $A$, the symbol $\|A\|_{\rm HS}$ denotes its Hilbert-Schmidt norm: $\|A\|_{\rm HS}=({\rm tr}(A^*A))^{\frac 12}$.
\end{itemize}
 
\bigskip

\section{Preliminaries}

Our proof of Theorem \ref{yap} uses complete integrability of equation \eqref{nls}. In that framework, \eqref{nls} can be solved by using the method of inverse scattering which we discuss next following \cite{FTbook}.

\subsection{The inverse scattering approach to NLS} Given a complex-valued function $q \in \Sch(\R)$, define the  differential operator
\begin{equation}\label{dir1}
L_{q} = i\sigma_3\frac{d}{d \xi} + i(q\sigma_- - \ov{q} \sigma_+),
\end{equation}
where we borrow notation  for constant matrices $\sigma_3$, $\sigma_{\pm}$  from \cite{FTbook}:
\begin{equation}\label{cmatrices}
\sigma_3 =  \begin{pmatrix}1 & 0 \\ 0 & -1\end{pmatrix}, 
\qquad
\sigma_+ = \begin{pmatrix}0&1\\0&0\end{pmatrix},
\qquad
\sigma_- = \begin{pmatrix}0&0\\1&0\end{pmatrix},
\qquad
\sigma = \sigma_- + \sigma_+=\begin{pmatrix}0&1\\1&0\end{pmatrix}.
\end{equation}
The expression $L_q$ is one of the forms in which the Dirac operator can be written. In {S}ection 3, we will introduce another form and will show how the two are related.
Let us also define 
$$
E(\xi, \lambda) = e^{\frac{\lambda}{2i}\xi\sigma_3} = \begin{pmatrix}e^{\frac{\lambda}{2i}\xi} & 0 \\ 0 & e^{-\frac{\lambda}{2i}\xi}\end{pmatrix}, 
$$
as in \cite{FTbook}. In the free case when $q=0$, the matrix-function $E$ solves $L_0E={\frac{\lambda}{2}} E,\, E(0,\lambda)=I$. 
Since $q \in \Sch(\R)$, it decays at infinity fast and therefore one can find two solutions
 $T_{\pm} = T_{\pm}(\xi, \lambda)$ such that 
\begin{equation}\label{eq1lol}
L_q T_{\pm} = \frac{\lambda}{2}T_{\pm}, \qquad T_{\pm} = E(\xi,\lambda) + o(1), \qquad \xi \to \pm\infty,
\end{equation}
for every $\lambda\in \R$.
These solutions are called the {\it Jost solutions} for $L_q$. Since both $T_{+}$ and $T_{-}$ solve the same ODE, they must satisfy
\begin{equation}\label{fix1}
T_{-}(\xi,\lambda) =  T_{+}(\xi,\lambda)T(\lambda), \qquad \xi \in \R, \qquad \lambda\in \R,
\end{equation}
where the matrix $T = T(\lambda)$ does not depend on $\xi \in \R$.  One can show that it has the form
\begin{equation}\label{eq4}
T(\lambda) = 
\begin{pmatrix}
a(\lambda) & \ov{b(\lambda)} \\
b(\lambda) & \ov{a(\lambda)}
\end{pmatrix}, \qquad \det T = |a|^2 - |b|^2 = 1.
\end{equation}
The matrix $T$ is called the {\it reduced transition matrix} for $L_q$, and the ratio $\rc_q = b/a$ is called the {\it reflection coefficient} for $L_q$. One can obtain $T$ in a different way: let $Z_q = Z_q(\xi, \lambda)$, $\xi \in \R$, $\lambda \in \C$ be the fundamental matrix  for $L_q$, that is, 
\begin{equation}\label{eq2}
L_q Z_q = \frac{\lambda}{2}Z_q, \qquad Z_q(0, \lambda) = \idm. 
\end{equation}
Then, we have
$Z_q(\xi,\lambda) = T_\pm(\xi,\lambda)T^{-1}_\pm(0,\lambda)$
and the pointwise limits
\begin{equation}\label{eq3}
T^{-1}_\pm(0,\lambda) = \lim_{\xi \to \pm\infty}E^{-1}(\xi,\lambda)Z_q(\xi,\lambda)
\end{equation}
exist for every $\lambda \in \R$. Moreover, we have $T(\lambda) = T_+^{-1}(0,\lambda)T_-(0,\lambda)$ on $\R$. 

\medskip

The coefficients $a,b$, and $\rc_q$ were defined for $\lambda\in \R$ and they satisfy $|a|^2=1+|b|^2$, $1-|\rc_q|^2=|a|^{-2}$ for these $\lambda$. However, one can show that $a(\lambda)$ is the boundary value of the outer function defined in ${\C}_+ = \{z \in \C:\; \Im z > 0\}$ by the formula (see (6.22) in \cite{FTbook})
\[
a(z)=\exp \left(\frac{1}{\pi i}\int_{\R}\frac{1}{\lambda-z} \log|a(\lambda)|\,d\lambda\right), \quad z\in \C_+,
\]
which, in view of identity $1-|\rc_q|^2=|a|^{-2}$ on $\R$, can be written as
\begin{equation}\label{ar}
a(z)=\exp \left(-\frac{1}{2\pi i}\int_{\R}\frac{1}{\lambda-z} \log(1-|\rc_q(\lambda)|^2)d\lambda\right)\,.
\end{equation}
That shows, in particular, that $b$ defines both $a$ and $\rc_q$, and $\rc_q$ defines $a$ and $b$.


\medskip

 The map $q\mapsto \rc_q$ is called the direct scattering transform and its inverse is called the inverse scattering transform. These maps are well-studied when $q\in \Sch(\R)$. In particular, we have the following result (see \cite{FTbook} for the proof).
\begin{Thm}\label{Thm1}
The map $q \mapsto \rc_q$ is a bijection from $\Sch(\R)$ onto the set of complex-valued functions $\{\rc\in\Sch(\R), \|\rc\|_{L^\infty(\R)}<1\}$.
\end{Thm}

The scattering transform has some symmetries: 

\begin{Lem}\label{sym} If $q\in \Sch(\R)$ and $\lambda \in \R$, then
\begin{eqnarray*}
 \text{\rm (dilation):} & \qquad \rc_{\alpha q(\alpha \xi)}(\lambda)=\rc_{q(\xi)}(\alpha^{-1}\lambda), \quad \alpha>0\,,
\\
\text{\rm (conjugation):}& \qquad \rc_{ \overline{q}(\xi)}(\lambda)={\ov{\rc_{q(\xi)}(-\lambda)}}\,,
\\
\text{\rm (translation):}& \qquad \rc_{q(\xi-\ell)}(\lambda)=\rc_{q(\xi)}(\lambda)e^{-i\lambda \ell}, \quad \ell\in \R\,,
\\
\text{\rm (modulation):}&  \qquad \rc_{e^{-i\beta \xi}q(\xi)}(\lambda)=\rc_{q(\xi)}(\lambda+\beta), \quad \beta\in \R\,.
\\
\text{\rm (rotation):}&  \qquad \rc_{\mu q(\xi)}(\lambda)={\mu}\rc_{q(\xi)}(\lambda), \quad \mu\in \C, \; |\mu|=1\,.
\end{eqnarray*}
\end{Lem}
\beginpf
Indeed, the direct substitution into \eqref{eq1lol} shows that  if $T_{\pm}(\xi,\lambda)$ are Jost solutions for $q(\xi)$, then 
\begin{itemize}
\item[$(a)$] $T_{\pm}(\alpha\xi,\alpha^{-1}\lambda)$ are the Jost solutions for $\alpha q(\alpha\xi)$,
\item[$(b)$] $\overline{T}_{\pm}(\xi,{-\lambda})$ are the Jost solutions for $\overline{q(\xi)}$,
\item[$(c)$] $T_{\pm}(\xi-\ell, {\lambda})E(\ell,\lambda)$ are the Jost solutions for $q(\xi-\ell)$,
\item[$(d)$] $E(-\xi,\beta)T_{\pm}(\xi,\lambda+\beta)$ are the Jost solutions for $e^{-i\beta \xi}q(\xi)$,
\item[$(e)$] $\left( \begin{smallmatrix}1 &0\\0&\mu\end{smallmatrix} \right)T_{\pm}(\xi,\lambda)\left( \begin{smallmatrix}1 &0\\0&\overline{\mu}\end{smallmatrix} \right)$ are the Jost solutions for $\mu q(\xi)$, $|\mu| = 1$.
\end{itemize}
Now, it is left to use the formula \eqref{fix1} which defines $T$. {A computation using \eqref{eq4} shows how $a$ and $b$ change under symmetries $(a)$--$(e)$. For example, the translation does not change $a$ and it multiplies $b$ by $e^{-i\lambda l}$. The modulation $e^{-i\beta \xi}q(\xi)$, however, gives $a_{e^{-i\beta \xi}q(\xi)}(\lambda)=a_{q(\xi)}(\lambda+\beta)$. Then, the claim follows from the definition of the reflection coefficient $\rc_q = b/a$.} 
\qed

\smallskip

The next result  {(see formula $(7.5)$ in \cite{FTbook})}, along with the previous theorem, shows how the inverse scattering transform can be used to solve \eqref{nls}.
\begin{Thm}\label{teo}
Let $q_0 \in \Sch(\R)$ and let $\rc_{q_0} = \rc_{q_0}(\lambda)$ be the reflection coefficient of $L_{q_0}$. Define the family 
\begin{equation}\label{con1}
\rc(\lambda, t) = e^{-i\lambda^2 t}\rc_{q_0}(\lambda), \qquad \lambda \in \R, \quad t \in \R.
\end{equation} 
For each $t \in \R$, let $q = q(\xi,t)$ be the potential in the previous theorem generated by $\rc(\lambda, t)$. Then, $q = q(\xi,t)$ is the unique classical solution of \eqref{nls} with the initial datum $q_0$. Moreover, for every $t \in \R$, the function $\xi \mapsto q(\xi, t)$ lies in $\Sch(\R)$.
\end{Thm}
The solutions to NLS equation
\begin{equation}\label{nls1}
i \frac{\partial q}{\partial t} = -\frac{\partial^2 q}{\partial \xi^2} + 2|q|^2q
\end{equation}
behave in an explicit way  {under some transformations}.  Specifically, we have  \medskip
\begin{itemize}
\item[$(a)$] Dilation: if $q(\xi,t)$ solves \eqref{nls1}, then $\alpha q(\alpha\xi,\alpha^2t)$ solves \eqref{nls1} for every $\alpha\neq 0$.
\item[$(b)$] Time reversal: if $q(\xi,t)$ solves \eqref{nls1}, then $\overline{q}(\xi,-t)$ solves \eqref{nls1}. In particular, if $q_0$ is real-valued, then $q(\xi,t)=q(\xi,-t)$.
\item[$(c)$] Translation: if $q(\xi,t)$ solves \eqref{nls1}, then $q(\xi-\ell,t)$ solves \eqref{nls1} for every $\ell\in \R$.
\item[$(d)$] Modulation or Galilean symmetry: if $q(\xi,t)$ solves \eqref{nls1}, then $e^{iv\xi-iv^2t}q(\xi-2vt,t)$ solves \eqref{nls1} for every $v\in \R$.
\item[$(e)$] Rotation: if $q(\xi,t)$ solves \eqref{nls1}, then $\mu q(\xi,t)$ solves \eqref{nls1} for every $\mu\in \C, |\mu|=1$.

\end{itemize}
These properties can be checked by direct calculation (see, e.g., formula (1.19) in \cite{kenig} for $(d)$) 
and a simple inspection shows that the bound \eqref{eq19} is consistent with all these transformations. The statements of Theorem \ref{teo} and  Lemma \ref{sym} are consistent with these symmetries as well.
\smallskip

Now, we can explain the idea behind the proof of the Theorem \ref{t5}.

\medskip

\noindent {\bf The idea of the proof for Theorem \ref{t5}.} One can proceed as follows. First, we assume that $q_0\in \Sch(\R)$ and notice that conservation of  {$|\rc(\lambda,t)|$,} $\lambda \in \R$, guaranteed by \eqref{con1}, yields that $\log |a(i,t)|$ is conserved, where $a(z,t)$ is defined for $z\in \C_+$ by \eqref{ar}. Separately,  for every Dirac  {operator $L_{q}$} with $q\in L^2(\R)$, we show that $\log|a(i)|$ is equivalent to some explicit quantity  {$\widetilde\K_{Q}$} that involves $q$. That quantity was introduced and studied in \cite{BD2017,BD2019,BD2022}: it resembles the matrix Muckenhoupt $A_2(\R)$  condition and it is equivalent to $H^{-1}(\R)$ norm of $q$ provided that $\|q\|_{L^2(\R)}$ is under control, e.g.,  
$
\|q\|_{L^2(\R)}<C
$
 with some fixed $C$. Putting things together, we see that Sobolev $H^{-1}(\R)$ norm  {of $q(\cdot, t)$} does not change much in time provided that the bound  {$
\|q(\cdot, t)\|_{L^2(\R)}<C
$} holds. Since  {$\|q(\cdot, t)\|_{L^2(\R)} = \|q_0\|_{L^2(\R)}$ is  time-invariant, we arrive to the statement of Theorem \ref{t5} for $q_0\in \Sch(\R)$ and $s=-1$. For $s=0$, the claim of Theorem \ref{t5} is trivial. The intermediate range of $s\in (-1,0)$ is handled by interpolation using Galilean invariance of NLS. The general case when $q_0\in L^2(\R)$ follows by a density argument if one uses the stability of $L^2$-solutions guaranteed by Theorem~\ref{yap}.}

\medskip

There are other methods that use conserved quantities that agree with negative Sobolev norms. The paper \cite{Killip18} uses a representation of $\log |a(i)|$ through a perturbation determinant. Then, the analysis of the perturbation series allows the authors of \cite{Killip18} to obtain estimates similar to \eqref{koch}. It is conceivable that this approach can provide results along the same lines as Theorem \ref{t5}.

\bigskip

To focus on the Dirac operator with $q\in L^2(\R)$, we first consider this operator on half-line $\R_+$ in connection to Krein systems that were introduced in \cite{Kr55}.
\subsection{Operator $L_q$ and Krein system} Let $A: \R_+ \to \C$ be a function on the positive half-line $\R_+ = [0, +\infty)$ such that
$$
\int_{0}^{r}|A(\xi)|\,d\xi < \infty, 
$$
for every $r \ge 0$. Recall that we denote the set of such functions by $L^{1}_{\loc}(\R_+)$. 
The Krein system (see the formula $(4.52)$ in \cite{Den06}) with the coefficient $A$ has the form  
\begin{equation}\label{eq7}
\begin{cases}
P'(\xi,\lambda) = i\lambda P(\xi,\lambda) - \overline{A(\xi)}P_*(\xi,\lambda), & P(0, \lambda) = 1, \\
P'_{*}(\xi,\lambda) = - A(\xi)P(\xi,\lambda), & P_*(0, \lambda) = 1,
\end{cases}
\end{equation}
where the derivative is taken with respect to $\xi \in \R_+$ and $\lambda\in \C$. Let also 
\begin{equation}\label{eq8}
\begin{cases}
\widehat{P}'(\xi,\lambda) = i\lambda \widehat{P}(\xi,\lambda) + \overline{A(\xi)}\widehat{P}_*(\xi,\lambda), & \widehat{P}(0, \lambda) = 1, \\
\widehat{P}'_{*}(\xi,\lambda) =  A(\xi)\widehat{P}(\xi,\lambda), & \widehat{P}_*(0, \lambda) = 1,
\end{cases}
\end{equation}
denote the so-called dual Krein system (see Corollary 5.7 in \cite{Den06}). Set 
\begin{equation}\label{eq6}
Y(\xi, \lambda) = e^{-i\lambda \xi} \begin{pmatrix}
P(2\xi, \lambda) & i\widehat{P}(2\xi, \lambda) \\
P_{*}(2\xi,\lambda) & -i\widehat{P}_{*}(2\xi,\lambda) 
\end{pmatrix}.
\end{equation}
The matrix-function  $Z_q$, which was defined in \eqref{eq2} for $q\in \Sch(\R)$, makes sense if we assume that $q\in L^1_{\rm loc}(\R)$. 
In the next lemma, we relate $Y$ to $Z_q$.
\begin{Lem}\label{l3}
Let $q \in L^1_{\loc}(\R)$,  $A(2\xi) = -\ov{q(\xi)}/2$ on $\R_+$, and  $Y$ be the corresponding matrix-valued function defined by \eqref{eq6}. Then, $Z_{q}(\xi, 2\lambda) = \sigma Y(\xi,\lambda)Y^{-1}(0,\lambda)\sigma$ for $\xi \ge 0$ and $\lambda \in \C$. 
\end{Lem}
\beginpf The proof is a computation.  We have
\begin{align*}
L_{\bar q} Y(\xi, \lambda) 
= & \lambda \sigma_{3} Y(\xi, \lambda) + i\sigma_3 e^{-i\lambda \xi}\frac{d}{d\xi}\begin{pmatrix}
P(2\xi, \lambda) & i\widehat{P}(2\xi, \lambda) \\
P_{*}(2\xi,\lambda) & -i\widehat{P}_{*}(2\xi,\lambda) 
\end{pmatrix} + i(\bar q\sigma_- - q \sigma_+)Y(\xi,\lambda),\\
= & 2i\sigma_3 e^{-i\lambda \xi}\begin{pmatrix}
i\lambda P(2\xi,\lambda) - \overline{A(2\xi)}P_*(2\xi,\lambda) & -\lambda \widehat{P}(2\xi,\lambda) + i\overline{A(2\xi)}\widehat{P}_*(2\xi,\lambda) \\
- A(2\xi)P(2\xi,\lambda) & -iA(2\xi)\widehat{P}(2\xi,\lambda) 
\end{pmatrix} \\
&+ i(\bar q\sigma_- - q \sigma_+ - i \lambda \sigma_{3})Y(\xi,\lambda).
\end{align*}
Notice, that
\begin{align*}
i(\bar q\sigma_- - q \sigma_+ - i \lambda \sigma_{3})Y(\xi,\lambda) 
&=ie^{-i\lambda \xi}
\begin{pmatrix}
-i\lambda & -q\\
\bar q & i\lambda
\end{pmatrix}
\begin{pmatrix}
P(2\xi, \lambda) & i\widehat{P}(2\xi, \lambda) \\
P_{*}(2\xi,\lambda) & -i\widehat{P}_{*}(2\xi,\lambda) 
\end{pmatrix}\\
&=
ie^{-i\lambda \xi}
\begin{pmatrix}
-i\lambda P(2\xi, \lambda) -qP_{*}(2\xi,\lambda)  & \lambda\widehat{P}(2\xi, \lambda) + iq\widehat{P}_{*}(2\xi,\lambda) \\
\bar qP(2\xi, \lambda) +i\lambda P_{*}(2\xi,\lambda) & i\bar q \widehat{P}(2\xi,\lambda) +\lambda\widehat{P}_{*}(2\xi,\lambda)
\end{pmatrix}\,.
\end{align*}
Using relation $2A(2\xi)+ \bar q(\xi) = 0$, we obtain
$$
L_{\bar q} Y(\xi, \lambda) = i e^{-i\lambda \xi} 
\begin{pmatrix}
i\lambda P(2\xi, \lambda)  & -\lambda\widehat{P}(2\xi, \lambda) \\
i\lambda P_{*}(2\xi,\lambda) & \lambda\widehat{P}_{*}(2\xi,\lambda) 
\end{pmatrix} = -\lambda Y(\xi, \lambda).
$$
Since $\sigma \sigma_3 \sigma = -\sigma_3$ and $\sigma \sigma_\pm \sigma = \sigma_{\mp}$, one has $\sigma L_{\bar q} \sigma = - L_{q}$. Therefore, 
$$
L_{q}(\sigma Y(\xi,\lambda)\sigma) = \lambda (\sigma Y(\xi,\lambda)\sigma).
$$
It follows that matrix-valued functions $Z_{q}(\xi, 2\lambda)$ and $\sigma Y(\xi,\lambda)Y^{-1}(0,\lambda)\sigma $ solve the same Cauchy problem and thus $Z_{q}(\xi, 2\lambda) = \sigma Y(\xi,\lambda)Y^{-1}(0,\lambda)\sigma$, as required. \qed 

\medskip

\begin{Lem}\label{l4}
Let $q \in L^1_{\loc}(\R)$, let $A(2\xi) = q(-\xi)/2$ on $\R_+$, and let $Y$ be the corresponding matrix-valued function defined by \eqref{eq6}. Then, $Z_{q}(-\xi, 2\lambda) = Y(\xi,\lambda) Y^{-1}(0,\lambda)$ for $\xi \ge 0$ and $\lambda \in \C$.
\end{Lem}
\beginpf  {Recall that matrices $\sigma_{3}$, $\sigma_{\pm}$, $\sigma$ are defined in \eqref{cmatrices}. Using  relations $\sigma \sigma_3 \sigma = -\sigma_3$ and $\sigma \sigma_\pm \sigma = \sigma_{\mp}$,} we see that $L_{\tilde q}\widetilde Z_{q} = \frac{\lambda}{2}\widetilde Z_{q}$, where $\tilde q(\xi) = -\ov{q(-\xi)}$ and $\widetilde Z_{q}(\xi, \lambda) = \sigma Z_{q}(-\xi, \lambda)\sigma$. Then, previous lemma applies to $\tilde q$, $Z_{\tilde q}(\xi, 2\lambda) = \widetilde Z_{q}(\xi, 2\lambda)$ and $A(2\xi) = -\ov{\tilde q(\xi)}/2 = q(-\xi)/2$. It gives $\widetilde Z_{q}(\xi, 2\lambda) = \sigma Y(\xi,\lambda)Y^{-1}(0,\lambda)\sigma$. Returning to $Z_q$, we get $Z_{q}(-\xi, 2\lambda) =Y(\xi,\lambda)Y^{-1}(0,\lambda)$. \qed

\medskip

\noindent Given $q \in L^2(\R)$, we define the continuous analogs of Wall polynomials (see \cite{KH01} and Section 7 in \cite{Den06}) by 
\begin{equation}\label{eq9}
\fA^{\pm} = \frac{P^{\pm}_* + \widehat{P}^{\pm}_*}{2}, \qquad \fA^{\pm}_* = \frac{P^{\pm} + \widehat{P}^{\pm}}{2}, \qquad 
\fB^{\pm} = \frac{P^{\pm}_* - \widehat{P}^{\pm}_*}{2}, \qquad \fB^{\pm}_* = \frac{P^{\pm}- \widehat{P}^{\pm}}{2},
\end{equation}
where $P^{\pm}$, $P^{\pm}_{*}$, $\widehat P^{\pm}$, $\widehat P^{\pm}_{*}$ are the solutions of systems \eqref{eq7}, \eqref{eq8} for the coefficient $A^+(\xi) = -\ov{q(\xi/2)}/2$ from Lemma \ref{l3} and the coefficient $A^{-}(\xi) = q(-\xi/2)/2$ from Lemma \ref{l4}, correspondingly. Functions $P^{\pm}$, $P^{\pm}_{*}$, $\widehat P^{\pm}$, $\widehat P^{\pm}_{*}$ are continuous analogs of polynomials orthogonal on the unit circle, they depend on two parameters: $\xi \in \R_+$ and $\lambda \in \C$ and they satisfy identities (see formula (4.32) in \cite{Den06}):
\begin{equation}\label{zam1}
P_{*}^{\pm}(\xi, \lambda) = e^{i\xi\lambda}\ov{P^{\pm}(\xi, \lambda)},\qquad \widehat P_{*}^{\pm}(\xi, \lambda) = e^{i\xi\lambda}\ov{\widehat P^{\pm}(\xi, \lambda)}
\end{equation}
for real $\lambda$.

\medskip

We will use the following result (see Lemma~2 in \cite{SA_GAFA} which contains a stronger statement). 
\begin{Thm}\label{t23}
Let $A \in L^2(\R_+)$, and let $P$, $P_{*}$ be the solutions of system \eqref{eq7} for the coefficient $A$. Then, the limit \begin{equation}
\Pi(\lambda) = \lim_{\xi \to +\infty}P_*(\xi, \lambda)\label{rest}
\end{equation}
exists for every $\lambda \in \C_+$. That function $\Pi$ is outer in $\C_+$. If $\lambda\in \R$, the convergence in \eqref{rest} holds  in the Lebesgue measure on $\R$ where $\Pi( \lambda)$ now denotes the non-tangential boundary value of $\Pi$.
\end{Thm}   
\noindent That theorem allows us to define 
\begin{equation}\label{eq10}
\fa^{\pm}(\lambda) = \lim_{\xi \to +\infty}  \fA^{\pm}(\xi,\lambda), 
\qquad 
\fb^{\pm}(\lambda) = \lim_{\xi \to +\infty}\fB^{\pm}(\xi,\lambda)
\end{equation}
for every $\lambda \in \C_+$ and for almost every $\lambda \in \R$.  Moreover,  Corollary 12.2 in \cite{Den06} gives
\begin{equation}\label{sa5}
|\fa^{\pm}(\lambda)|^2=1+|\fb^{\pm}(\lambda)|^2
\end{equation}
for a.e. $\lambda\in \R$.  {For every $\lambda\in \C_+$, we define
\[
a(\lambda)=\fa^+(\lambda) \fa^-(\lambda) - \fb^+(\lambda)\fb^-(\lambda)\,.
\]
\begin{Prop}\label{outer}
The function $a$ is outer in $\C_+$.
\end{Prop}
\beginpf
We can write
\[
a=\fa^+ \fa^-(1 - s^+s^-), \quad s^{\pm}=\frac{\fb^\pm}{\fa^\pm}\,.
\]
It is known that $\fa^\pm$ are outer (see the formulas (12.9) and (12.29) in \cite{Den06}) and that $s^\pm$ satisfy $|s^\pm|<1$ in $\C_+$. The function $1 - s^+s^-$ has a positive real part in $\C_+$ and so is an outer function. That shows that $a$ is a product of three outer functions and hence it is outer itself.
\qed
}
\medskip

\begin{Prop}\label{p21}
Let $q \in L^2(\R)$ and let $Z_q$ be defined by \eqref{eq2}. Then, the limits in \eqref{eq3} exist in the Lebesgue measure on $\R$.
The matrix $T(\lambda) = T_+^{-1}(0,\lambda)T_-(0,\lambda)$ has the form \eqref{eq4} where
\begin{equation}\label{eqab}
a = \fa^+ \fa^- - \fb^+\fb^-, \qquad b = \fa^{-}\ov{\fb^+} - \fb^{-}\ov{\fa^+}\,,
\end{equation}
and $\fa^{\pm}$, $\fb^{\pm}$ are defined Lebesgue almost everywhere on $\R$ by the convergence in \eqref{eq10} in measure.
\end{Prop}
\beginpf If $q \in L^2(\R)$, the fundamental matrix $Z_q$ and the continuous Wall polynomials \eqref{eq9} are related by the formula
\begin{eqnarray}\label{refe}
Z_{q}(\xi, 2\lambda) 
= 
\begin{cases}
e^{-i\lambda\xi}
\begin{pmatrix}
\fA^+(2\xi, \lambda) & \fB^+(2\xi,\lambda)\\
\fB^+_*(2\xi, \lambda) & \fA^+_*(2\xi,\lambda)
\end{pmatrix}, & \xi \ge 0,\\
e^{i\lambda\xi}
\begin{pmatrix}
\fA^-_*(-2\xi,\lambda) & \fB^-_*(-2\xi, \lambda)\\
\fB^-(-2\xi,\lambda) & \fA^-(-2\xi, \lambda)
\end{pmatrix}, & \xi < 0.
\end{cases}
\end{eqnarray}
Indeed, it is enough to use Lemma \ref{l3}, Lemma \ref{l4}, and the fact that $Y^{-1}(0, \lambda) = \frac{1}{2}\left(\begin{smallmatrix}1 & 1\\-i &i\end{smallmatrix}\right)$. Our next step is to prove that the limit
\begin{equation}
T_+^{-1}(0,2\lambda) = \lim_{\xi \to +\infty}E^{-1}(\xi,2\lambda)Z_q(\xi,2\lambda)
\end{equation}
exists in Lebesgue measure when $\lambda\in \R$. From \eqref{zam1}, we obtain
\begin{align*}
E^{-1}(\xi,2\lambda)Z_q(\xi,2\lambda) 
= \begin{pmatrix}1 & 0\\ 0& e^{-2i\lambda\xi}\end{pmatrix}
\begin{pmatrix}
\fA^+(2\xi, \lambda) & \fB^+(2\xi,\lambda)\\
\fB^+_*(2\xi, \lambda) & \fA^+_*(2\xi,\lambda)
\end{pmatrix}
=
\begin{pmatrix}
\fA^+(2\xi, \lambda) & \fB^+(2\xi,\lambda)\\
\ov{\fB^+}(2\xi, \lambda) & \ov{\fA^+}(2\xi,\lambda)
\end{pmatrix},
\end{align*}
for every $\xi \ge 0$ and $\lambda \in \R$.  Similarly, 
\begin{align*}
E^{-1}(-\xi,2\lambda)Z_q(-\xi,2\lambda) 
= \begin{pmatrix}e^{-2i\lambda\xi} & 0\\ 0& 1\end{pmatrix}
\begin{pmatrix}
\fA^-_*(2\xi,\lambda) & \fB^-_*(2\xi, \lambda)\\
\fB^-(2\xi,\lambda) & \fA^-(2\xi, \lambda)
\end{pmatrix}
=
\begin{pmatrix}
\ov{\fA^-}(2\xi, \lambda) & \ov{\fB^-}(2\xi,\lambda)\\
\fB^-(2\xi,\lambda) & \fA^-(2\xi, \lambda)
\end{pmatrix}.
\end{align*}
Hence, the limits 
\begin{equation}
T^{-1}_{\pm}(0,2\lambda)
= \lim_{\xi \to \pm\infty}E^{-1}(\xi,2\lambda)Z_q(\xi,2\lambda)
\end{equation}
exist in Lebesgue measure on $\R$  {by Theorem \ref{t23}}. Moreover, 
\begin{align*}
T(2\lambda) 
&= T_+^{-1}(0,2\lambda)T_-(0,2\lambda) 
= 
\begin{pmatrix}
\fa^+(\lambda) & \fb^{+}(\lambda) \\
\ov{\fb^{+}(\lambda)} & \ov{\fa^{+} (\lambda)}
\end{pmatrix} \begin{pmatrix}
\ov{\fa^{-}(\lambda)} & \ov{\fb^{-}(\lambda)} \\
\fb^-(\lambda) & \fa^{-}(\lambda)
\end{pmatrix}^{-1}\\
&\stackrel{\eqref{sa5}}{=}
\begin{pmatrix}
\fa^+(\lambda) & \fb^{+}(\lambda) \\
\ov{\fb^{+}(\lambda)} & \ov{\fa^{+} (\lambda)}
\end{pmatrix} 
\begin{pmatrix}
\fa^{-}(\lambda) & -\ov{\fb^-(\lambda)}\\
-\fb^{-}(\lambda) & \ov{\fa^{-}(\lambda)}
\end{pmatrix}
=
\begin{pmatrix}
a(\lambda) & \ov{b(\lambda)} \\
b(\lambda) & \ov{a(\lambda)}
\end{pmatrix}
\end{align*}
and the proposition follows. \qed

\medskip

We end this section with a few remarks on reflection coefficients of potentials in $L^2(\R)$. We have $|a|^2 - |b|^2 = 1$ almost everywhere on $\R$ due to the fact that $\det T_{\pm}(0, \lambda) = 1$ almost everywhere on $\R$. That can also be established directly using \eqref{sa5}. Proposition \ref{p21} then allows to define the reflection coefficient $\rc_{q} = b/a$ for every $q \in L^2(\R)$. The Lemma \ref{sym} holds for $\rc_{q}$ in that case as well. However, not all results about scattering transform can be generalized from the case $q\in \Sch(\R)$ to $q\in L^2(\R)$. For example, scattering transform is injective on $\Sch(\R)$ by Theorem~\ref{Thm1}, but it is not longer so when extended to $L^2(\R)$ (see Example \ref{ex1} in Appendix).\bigskip

\section{Another form of Dirac operator, $q\in L^2(\R)$, and  the entropy function.}\label{s3}

Suppose $q\in L^2(\R)$. The alternative to $L_q$ form of writing Dirac operator on the line is given by an expression
\begin{equation}\label{eq18}
\Di_{Q}: X \mapsto JX' + QX, \qquad Q = \begin{pmatrix}-\Im q & -\Re q \\ -\Re q & \Im q\end{pmatrix}\,.
\end{equation} 
$\Di_Q$  is densely defined self-adjoint operator on the Hilbert space $L^2(\R, \C^2)$ of functions $X: \R \to \C^2$ such that $\|X\|_{L^2(\R, \C^2)}^{2} = \int_{\R}\|X(\xi)\|^2_{\C^2}\,d\xi$ is finite. $\Di_Q$ and $L_{q}$ defined in \eqref{dir1} are related by a simple formula: 
\[
\Di_Q=\Sigma L_{q} \Sigma^{-1}, \quad \Sigma=\frac{1}{\sqrt 2}\left( \begin{matrix} 1 &1\\ -i &i \end{matrix} \right), \quad \Sigma^{-1}=\frac{1}{\sqrt 2}\left( \begin{matrix} 1 &i\\ 1 &-i \end{matrix} \right)\,.
\]
One way to study $\Di_Q$ is to  focus on Dirac operators on half-line $\R_+$ first.
Given  $q \in L^2(\R_+)$, we define $\Di^+_{Q}$ on $L^2(\R_+, \C^2)$ by 
\begin{equation}\label{eq14}
\Di^{+}_{Q}: X \mapsto JX' + QX, \qquad Q = \begin{pmatrix}-\Im q & -\Re q \\ -\Re q & \Im q\end{pmatrix}
\end{equation}
on the dense subset of absolutely continuous functions $X \in L^2(\R_+, \C^2)$ such that $\Di^{+}_{Q}X \in L^2(\R_+, \C^2)$, $X(0) = \left(\begin{smallmatrix}*\\ 0\end{smallmatrix}\right)$. We will call $\Di^+_Q$ the Dirac operator defined on the positive half-line with boundary conditions $X(0) = \left(\begin{smallmatrix}*\\ 0\end{smallmatrix}\right)$ or simply the half-line Dirac operator. Set $A(\xi) = -\ov{q(\xi/2)}/2$ for $\xi \in \R_+$, and let $P(\xi,\lambda)$, $P_*(\xi,\lambda)$ be the solutions of Krein system \eqref{eq7} generated by $A$. The Krein system with coefficient $A$ and Dirac equation \eqref{eq14} are related as follows (see the proof of Lemma \ref{l6} in Appendix): if $N_Q$ solves the Cauchy problem $JN'_Q(\xi,\lambda) + Q(\xi) N_Q(\xi,\lambda) = \lambda N_Q(\xi,\lambda)$, $N_Q(0,\lambda) = \idm$, then 
\[
N_Q(\xi,\lambda)= \frac{e^{-i\lambda \xi}}{2} \left(\begin{matrix} 1&1\\i&-i\end{matrix}\right) \left(\begin{matrix}  \fA_*^+(2\xi,\lambda)&\fB_*^+(2\xi,\lambda)\\
\fB^+(2\xi,\lambda)&\fA^+(2\xi,\lambda) \end{matrix}\right)\left(\begin{matrix} 1&-i\\1&i\end{matrix}\right)\,,
\]
where the continuous Wall polynomials $\fA^+,\fB^+,\fA_*^+,\fB^+_*$ were defined in \eqref{eq9}. The Weyl function of the operator $\Di^{+}_{Q}$ coincides (see Lemma \ref{l6} in Appendix) with 
\begin{equation}\label{eq15}
m_{Q}(z) = \lim_{\xi \to +\infty} i\frac{\widehat P_{*}(\xi,z)}{P_{*}(\xi,z)}, \qquad z \in \C_+.
\end{equation}
It is known (see Theorem 7.3 in \cite{Den06}) that the limit above exists for every $z \in \C_+$ and defines an analytic function of Herglotz-Nevanlinna class in $\C_+$. The latter means that $m_{Q}(\C_+) \subset \C_+$. In the next theorem, $\Im m_Q(\lambda)$ denotes the nontangential boundary value on $\R$ which exists Lebesgue almost everywhere. It is understood as a nonnegative function $g = \Im m$ on $\R$ and it satisfies $g/(1+\lambda^2) \in L^1(\R)$.
\begin{Thm}\label{t15}
Let $q \in L^2(\R_+)$ and let $Q$, $\Di^{+}_{Q}$, $m_{Q}$ be defined by \eqref{eq14}, \eqref{eq15}. Denote by $N_Q$ the solution of the Cauchy problem $JN'_Q(\xi) + Q(\xi) N_Q(\xi) = 0$, $N_Q(0) = \idm$, and set $\Hh_{Q} = N_Q^*N_Q$. Define also 
\begin{align}
\K^{+}_{Q}
&= \log \Im m_{Q}(i) - \frac{1}{\pi}\int_{\R}\log \Im m_{Q}(\lambda)\frac{d\lambda}{\lambda^2 + 1},\label{eq16}\\
\widetilde \K^{+}_{Q} 
&= \sum_{k =0 }^{+\infty}\left(\det\int_{k}^{k+2}\Hh_Q(\xi)\,d\xi - 4\right).  \label{faw}
\end{align}
Then, we have
\begin{equation}\label{eq17}
c_1 \K_{Q}^{+} \le \widetilde\K_{Q}^{+} \le c_2 e^{c_2 \K_{Q}^{+}} 
\end{equation}  
for some positive absolute constants $c_1$, $c_2$.
\end{Thm} 
\beginpf
Lemma \ref{l6} in Appendix shows that $m_Q$ coincides with the Weyl function for the canonical system with Hamiltonian $\Hh_Q$.
Then, the bounds in \eqref{eq17} follow from the Theorem 1.2 in \cite{BD2019} (see also Corollary~1.4 in \cite{BD2019}).\qed
\medskip


The quantity $\K^+_Q$ will be called {\it the entropy} of the Dirac operator on {$\R_+$}. We now turn to \eqref{eq18} to define the entropy for the Dirac operator on the whole line. Take $q \in L^2(\R)$ and let
$A^+(\xi) = -\ov{q(\xi/2)}/2$ and $A^{-}(\xi) = q(-\xi/2)/2$, $\xi \in \R_+$ be the coefficients of Krein systems associated to restrictions of $q$ to the half-lines $\R_+$ and $\R_-$. As in \eqref {eq15}, the half-line Weyl functions $m_{\pm}$ satisfy 
\begin{equation}\label{eq23}
m_{\pm}(z) = \lim_{\xi \to +\infty} i\frac{\widehat P_{*}^{\pm}(\xi,z)}{P_*^{\pm}(\xi,z)}, \qquad z \in \C_+.
\end{equation}
These Weyl functions $m_\pm$ can be used to construct the spectral representation for the Dirac operator.
Let
\begin{equation}\label{eq22}
m(z) = -\frac{1}{m_+(z) + m_-(z)}
\begin{pmatrix}
-2m_+(z)m_-(z) & m_+(z) - m_-(z)\\
m_+(z) - m_-(z) & 2
\end{pmatrix}, \qquad z \in \C_+.
\end{equation}
Using $\Im m_{\pm}(z) > 0$, one can show that $\Im m(z)$ is a positive definite matrix for $z \in \C_+$. In other words, $m$ is the matrix-valued Herglotz function. Therefore, there exists a unique matrix-valued measure $\rho$ taking Borel subsets of $\R$ into $2 \times 2$ nonnegative matrices such that
$$
m(z) 
= \alpha + \beta z + \frac{1}{\pi}\int_{\R}\left(\frac{1}{\lambda-z}
- \frac{\lambda}{\lambda^2+1}\right)\,d\rho(\lambda), \qquad z \in \C_+, 
$$
where $\alpha$, $\beta$ are constant $2 \times 2$ real matrices, $\beta \ge 0$.  The importance of  $\rho$ becomes clear when we recall the spectral decomposition for $\Di_Q$. Specifically,  let $N_Q(\xi, z)$ be the solution of the Cauchy problem
\begin{equation}\label{eq20}
J\frac{\partial}{\partial \xi}N_Q(\xi, z) + Q(\xi)N_Q(\xi, z) = z {N}_{Q}(\xi, z), \quad \quad N_Q(0, z) = \idm, \quad z \in \C, \quad \xi\in \R\,.
\end{equation}
Then, the mapping
\begin{equation}\label{wt}
\F_{\Di_Q}: X \mapsto \frac{1}{\sqrt{\pi}}\int_{\R}N_{Q}^*(\xi, \lambda)X(\xi)\,d\xi, \qquad \lambda \in \R,
\end{equation}
initially defined on the set of compactly supported smooth functions $X: \R \to \C^2$, extends (see Appendix) to the unitary operator between the Hilbert spaces $L^2(\R, \C^2)$ and $L^2(\rho)$, 
\begin{align*}
L^2(\rho) 
&= \Bigl\{Y: \R \to \C^2: \; \|Y\|^2_{L^2(\rho)} = \int_{\R}Y^*(\lambda)\,d\rho(\lambda)\,Y(\lambda) < \infty\Bigr\}.
\end{align*}
Moreover, $\Di_{Q}$ is unitary equivalent to the operator of multiplication  by the independent variable in $L^2(\rho)$ and the unitary equivalence is given by the operator $\F_{Q}$.  {In fact, these properties of $\rho$  will not be used in the paper, we mention them  only to motivate the following definition.}   Let us define the {\it entropy function} $\K_Q(z)$ by
\begin{equation}\label{cot}
\K_{Q}(z) = -\frac{1}{\pi}\int_{\R}\log (\det \rho_{ac}(\lambda)) \frac{\Im z }{|\lambda - z|^2}\, d\lambda, \qquad z \in \C_+,
\end{equation}
where $\rho_{ac}$ denotes the absolutely continuous part of the spectral measure $\rho$ and it satisfies $\rho_{ac}(\lambda)=\lim_{\eps\to 0,\eps>0} \Im m(\lambda+i\eps)$ for a.e. $\lambda\in \R$. The quantity $\K_Q$ will play a crucial role in our considerations. We first relate it to the coefficient $a$ of the reduced transition matrix $T$ which was introduced in  Proposition~\ref{p21}.
\begin{Lem}\label{l2}
We have  {$\det \rho_{ac}(\lambda) = |a(\lambda)|^{-2}$} for almost all $\lambda \in \R$. In particular, $\K_{Q}(z) = 2\log |a(z)|$ for all $z \in \C_+$. 
\end{Lem}
\beginpf From the definition (or see page 59 in \cite{Remlingb}), one has
\begin{equation}\label{sn1}
\det \Im m(z) = 4\frac{\Im m_+(z) \Im m_-(z)}{|m_+(z) + m_-(z)|^2}, \qquad z \in \C_+.
\end{equation}
Substituting expressions for 
$$
m_{\pm}(z) = \lim_{\xi \to +\infty} i\frac{\widehat P_{*}^{\pm}(\xi,z)}{P^{\pm}(\xi,z)} = i\frac{\fa^{\pm}(z) - \fb^{\pm}(z)}{\fa^{\pm}(z) + \fb^{\pm}(z)}, \qquad z \in \C_+,
$$ 
into \eqref{sn1},
we obtain 
\begin{align*}
\det \rho_{ac}(\lambda) 
&= \lim_{\eps \to +0} \det \Im m(\lambda + i\eps) \\
&= \lim_{\eps \to +0} \frac{(|\fa^+(\lambda + i\eps)|^2 - |\fb^+(\lambda+i\eps)|^2)(|\fa^-(\lambda + i\eps)|^2 - |\fb^-(\lambda+i\eps)|^2)}{|\fa^{+}(\lambda + i\eps)\fa^-(\lambda + i\eps) - \fb^{+}(\lambda + i\eps)\fb^{-}(\lambda + i\eps)|^2}\\
&=\frac{1}{|a(\lambda)|^2},
\end{align*}
for almost every $\lambda \in \R$ and  {the first claim of the lemma follows. Then, the second claim is immediate because $a$ is an outer function as we showed in Proposition \ref{outer}.\qed}

\medskip

 Consider again the half-line entropy functions
$$
\K_{Q}^{\pm}(z) = \log \Im m_{\pm}(z) - \frac{1}{\pi}\int_{\R}\log \Im m_{\pm}(\lambda)\frac{\Im z}{|\lambda - z|^2}\,d\lambda, \qquad z \in \C_+.
$$
We see that $\K_{Q}^{+}(i)$ coincides with the entropy \eqref{eq16} for the restriction of $Q$ to $\R_+$ (that explains why we use the same notation for the two objects), and $\K^{-}_{Q}(z) = \K^{+}_{Q_-}(z)$ for the potential  
$$
Q_{-}(\xi) = \begin{pmatrix}-\Im q(-\xi) & \Re q(-\xi) \\ \Re q(-\xi) & \Im q(-\xi) \end{pmatrix}, \qquad 	\xi \in \R_+.
$$ 
Our plan now is to relate $\K^{\pm}_{Q}(i)$ with $\K_{Q}(i)$ and then use the fact that the full line entropy $\K_{Q}(i)$ is conserved, see Lemma \ref{l2}. That will eventually lead to the proof of Theorem \ref{t5}.

\begin{Lem}\label{l5}
Let $q \in L^2(\R)$ and let $q_{\ell}(\xi) = q(\xi - \ell)$, where  $\ell \in \R$ and $\xi \in \R$. Denote by $Q_{\ell}$ the matrix-function in \eqref{eq18} corresponding to $q_{\ell}$. Then, $\K^{+}_{Q_{\ell}}(z) \to \K_{Q}(z)$, $\K^{-}_{Q_{\ell}}(z) \to 0$ as $\ell \to +\infty$ for every $z \in \C_+$.
\end{Lem}
\beginpf Take $z \in \C_+$. We have
$$
\K^{+}_{Q_{\ell}}(z) + \K^{-}_{Q_{\ell}}(z) 
= 
\log \Bigl(\Im m_{\ell, +}(z)\Im m_{\ell, -}(z)\Bigr) - \frac{1}{\pi}\int_{\R}\log \Bigl( \Im m_{\ell,+}(\lambda)\Im m_{\ell,-}(\lambda)\Bigr)\frac{\Im z}{|\lambda - z|^2}\,d\lambda,
$$
for the corresponding Weyl functions $m_{\ell, \pm}$. We also have 
$$
\log |m_{\ell,+}(z) + m_{\ell,-}(z)|^2 = \frac{1}{\pi}\int_{\R}\log |m_{\ell,+}(\lambda) + m_{\ell,-}(\lambda)|^2\frac{\Im z}{|\lambda - z|^2}\,d\lambda
$$
by the mean value theorem for harmonic functions. From \eqref{sn1}, it follows that
$$
\K^{+}_{Q_\ell}(z) + \K^{-}_{Q_\ell}(z) =  \log \left(4\frac{\Im m_{\ell,+}(z)\Im m_{\ell,-}(z)}{|m_{\ell,+}(z) + m_{\ell,-}(z)|^2}\right) + \K_{Q_{\ell}}(z). 
$$
Notice that $\K_{Q_{\ell}}(z)$ does not depend on $\ell \in \R$ because the coefficient $a$ in Lemma \ref{l2} for the potential $Q_\ell$ does not depend on $\ell$.
So, we only need to show that
$$
\K^{-}_{Q_\ell}(z) \to 0 \quad\,\,{\rm and}\,\,\quad \log \left(4\frac{\Im m_{\ell,+}(z)\Im m_{\ell,-}(z)}{|m_{\ell,+}(z) + m_{\ell,-}(z)|^2}\right) \to 0,
$$
when $\ell \to +\infty$ and $z\in \C_+$. The second relation follows from $m_{\ell,+}(z) \to i$, $m_{\ell,-}(z) \to i$, which hold because  $q_{\ell}$ tends to zero weakly in $L^2(\R)$ as $\ell \to +\infty$  and $\|q_{\ell}\|_{L^2(\R)}=\|q\|_{L^2(\R)}$ (see Lemma \ref{l7} in Appendix). Moreover, relation $m_{\ell,-}(z) \to i$ implies  that 
$\K^{-}_{Q_\ell}(z) \to 0$ if and only if
\begin{equation}\label{low1}
\frac{1}{\pi}\int_{\R}\log\Im m_{\ell,-}(\lambda)\frac{\Im z}{|\lambda - z|^2}\,d\lambda \to 0\,.
\end{equation}
In the rest of the proof, we will show \eqref{low1}.
Let $\fa_{\ell}^{-}$, $\fb_{\ell}^{-}$ be the limits of continuous Wall polynomials corresponding to $Q_{\ell}^-$. Consider $s_\ell^-=\fb_\ell^-/\fa_\ell^-$. The formula (12.57) in \cite{Den06} gives 
\[
s_\ell^-(z)=\frac{1+im_{\ell,-}(z)}{1-im_{\ell,-}(z)}\,, \qquad  {m_{\ell,-}(z) = i \frac{\fa_\ell^-(z) - \fb^-_\ell(z)}{\fa^-_\ell(z) + \fb_\ell^-(z)}}.
\]
It implies that
$\Im m_{\ell,-}(\lambda)=|\fa^{-}_{\ell}(\lambda) + \fb^{-}_{\ell}(\lambda)|^{-2}$ when $\lambda\in \R$ and that $s^-_\ell(z)\to 0$ when $\ell\to +\infty$ and $z\in \C_+$. Now, we can write
\begin{align*}
\frac{1}{\pi}\int_{\R}\log\Im m_{\ell,-}(\lambda)\frac{\Im z}{|\lambda - z|^2}\,d\lambda 
&= \frac{1}{\pi}\int_{\R}\log\left(\frac{1}{|\fa^{-}_{\ell}(\lambda) + \fb^{-}_{\ell}(\lambda)|^2}\right)\frac{\Im z}{|\lambda - z|^2}\,d\lambda \\
&= \log \frac{1}{|\fa^{-}_{\ell}(z) + \fb^{-}_{\ell}(z)|^2} = \log \frac{1}{|\fa_{\ell}^{-}(z)|^2} + \log \frac{1}{|1 + s_{\ell}^{-}(z)|^2}\,.
\end{align*}
 {So, it remains to show that $|\fa_{\ell}^{-}(z)|^2 \to 1$ as $\ell \to+\infty$. That holds because $\|q_{\ell,-}\|_{L^2(\R_+)} \to 0$  as $\ell \to+\infty$ and
$$
\|q_{\ell,-}\|_{L^2(\R_+)}^{2} = \frac{1}{\pi}\int_{\R}\log|\fa_{\ell}^{-}(\lambda)|^2 \,d\lambda \ge
\frac{\Im z}{\pi} \int_{\R}\log|\fa_{\ell}^{-}(\lambda)|^2\frac{\Im z}{|\lambda - z|^2} \,d\lambda = \Im z\cdot \log|\fa_{\ell}^{-}(z)|^2 \ge 0,
$$
where the first equality follows from $\|q_{\ell,-}\|_{L^2(\R_+)}^2=2\|A_{\ell,-}\|_{L^2(\R_+)}^2$ and the formula $(12.2)$ in \cite{Den06}. Thus, \eqref{low1} holds and we are done.} \qed

\medskip

As an immediate corollary of Theorem \ref{t15} and Lemma \ref{l5}, we get the following estimate.

\begin{Lem}\label{l88}
Let $q \in L^2(\R)$. Denote by $N_Q$ the solution of the Cauchy problem $JN_Q'(\xi) + Q(\xi) N_Q(\xi) = 0$, $N_Q(0) = \idm$, and set $\Hh_{Q} = N_Q^*N_Q$. Consider
\begin{equation}\label{faw1}
\K_{Q}= \K_{Q}(i), \qquad 
\widetilde\K_{Q} = \sum_{k  \in \Z}\left(\det\int_{k}^{k+2}\Hh_Q(\xi)\,d\xi - 4\right).
\end{equation}
Then, we have
\begin{equation}\label{ler}
c_1 \K_{Q} \le \widetilde\K_{Q} \le c_2 \K_{Q} e^{c_2 \K_{Q}}
\end{equation}
for some positive absolute constants $c_1$, $c_2$. 
\end{Lem}
\beginpf By Lemma \ref{l5}, we have $\K_{Q} = \lim_{\ell \to +\infty}\K^{+}_{Q_\ell}$. 
It remains to substitute $Q_{\ell}$ into the  estimate \eqref{eq17} and take the limit as $\ell \to +\infty$ for $\ell\in \Z$.  \qed\bigskip

\section{Proof of Theorem \ref{t5}}
The following result will play a crucial role in what follows. We postpone its proof to the next section.
\begin{Thm}\label{P2c1b}
Suppose  $q\in L^2(\R)$ and let $N_Q$ satisfy $JN_Q' + QN_Q = 0, N_Q(0)=I$, where  $Q = \left(\begin{smallmatrix}-\Im q & -\Re q \\ -\Re q & \Im q\end{smallmatrix}\right).$ Then, 
\begin{eqnarray}\label{eq41}
e^{-C_1R}\|q\|_{H^{-1}(\R)}^2\lesssim \widetilde \K_{Q}\lesssim  e^{C_2R}\|q\|^2_{H^{-1}(\R)}\,,
\end{eqnarray}
where $R=\|q\|_{L^2(\R)}$ and $C_1$, $C_2$ are two positive absolute constants.
\end{Thm}
\smallskip

\noindent{\bf Proof of Theorem \ref{t5} in the case $s = -1$.}  First, assume that $q_0 \in \Sch(\R)$  and let $q(\xi, t)$ be the solution of \eqref{nls} with the initial datum $q_0$.  We want to prove that  
\begin{equation}\label{m1}
C_1(1+\|q_0\|_{L^2(\R)})^{-2} \|q_0\|_{H^{-1}(\R)}\le \|q(\cdot, t)\|_{H^{-1}(\R)} \le C_2(1+\|q_0\|_{L^2(\R)})^{2} \|q_0\|_{H^{-1}(\R)}\,.
\end{equation} 
We have $\|q(\cdot,t)\|_{L^2(\R)}=\|q_0\|_{L^2(\R)}$ for all $t$, see formula $(4.33)$ in  \cite{FTbook}.   Let $a(z,t)$ denote the coefficient in the matrix \eqref{eq4} given by $q(\xi,t)$.
For each $t \in \R$, define $Q$ by \eqref{eq14}. Let $\widetilde \K_{Q}(t)$ be as in Lemma~\ref{l88} and $\K_Q(z,t)$ be defined by \eqref{cot}.  Formulas \eqref{ar} and \eqref{con1} show that $a({z},t)$ is constant in $t$ and Lemma~\ref{l2} says that $\K_Q(z,t)$ is constant in $t$ as well. The bound \eqref{ler} yields
\begin{equation}\label{eq32}
c_1\K_Q(i,0)\le \widetilde \K_Q(t)\le c_2\K_Q(i,0)e^{c_2\K_Q(i,0)}\,.
\end{equation}

 {
Assume first that  $R=\|q_0\|_{L^2(\R)}\le 1$. Taking $t=0$ in \eqref{eq32} and applying \eqref{eq41} to $q_0$, we get $\K_Q(i,0)\lesssim 1$ since $\|q_0\|_{H^{-1}(\R)}\le R\le 1$. Hence, in that case \eqref{eq32} can be written as
$
 \widetilde \K_Q(t)\sim\K_Q(i,0)\,.
$
By \eqref{P2c1b}, $\|q(\cdot, t)\|_{H^{-1}(\R)}^2\sim \widetilde \K_Q(t)$, and so $\|q(\cdot, t)\|_{H^{-1}(\R)}^2\sim \|q(\cdot, 0)\|_{H^{-1}(\R)}^2$.}\smallskip

If $R=\|q_0\|_{L^2(\R)}> 1$, we  use dilation. Consider $q_\alpha(\xi,t)=\alpha q(\alpha\xi, \alpha^2t)$ which solves the same equation and notice that 
$
\|q_\alpha\|_{L^2(\R)}=\alpha^{{\frac 12}}R\,.
$

Let $\alpha=\alpha_c=R^{-2}<1$ making $\|q_{\alpha_c}\|_{L^2(\R)}=1$. Then, for the Sobolev norm, we get
\begin{equation}\label{lw1}
\|q_\alpha(\cdot,t)\|_{H^{-1}(\R)}=\alpha^{{\frac 12}}\left(\int_\R \frac{1}{1+\alpha^2\eta^2}|(\cal{F} q)(\eta,{\alpha^2}t)|^2d\eta\right)^{{\frac 12}}\,.
\end{equation}
Since 
\begin{equation}\label{lop1}
\frac{1}{1+\eta^2}\le \frac{1}{1+\alpha_c^2\eta^2}\le \frac{1}{\alpha_c^2(1+\eta^2)},
\end{equation}
one has
\[
\alpha_c^{{\frac 12}}\|{q(\cdot, \alpha_c^2 t)}\|_{H^{-1}(\R)}{\le} \|q_{\alpha_c}(\cdot,{t})\|_{H^{-1}(\R)}{\le} \alpha_c^{-{\frac 12}}\|{q(\cdot, \alpha_c^2 t)}\|_{H^{-1}(\R)}\,.
\]
{In particular, at $t = 0$ we get
\[
\alpha_c^{{\frac 12}}\|q(\cdot, 0)\|_{H^{-1}(\R)}
\le \|q_{\alpha_c}(\cdot,0)\|_{H^{-1}(\R)}
\le \alpha_c^{-{\frac 12}}\|q(\cdot, 0)\|_{H^{-1}(\R)}\,.
\]
Since $\|q_{\alpha_c}(\cdot,0)\|_{L^2(\R)}=1$, one can apply the previous bounds to obtain 
\[
\|q_{\alpha_c}(\cdot,t)\|_{H^{-1}(\R)}\sim  \|q_{\alpha_c}(\cdot,0)\|_{H^{-1}(\R)}\,.
\]
Then,
\begin{align*}
\alpha_c^{{\frac 12}}\|q(\cdot, \alpha_c^2 t)\|_{H^{-1}(\R)} 
&\lesssim \|q_{\alpha_c}(\cdot,0)\|_{H^{-1}(\R)} \lesssim 
\alpha_c^{-{\frac 12}}\|q(\cdot, 0)\|_{H^{-1}(\R)},\\
\alpha_c^{-{\frac 12}}\|q(\cdot, \alpha_c^2 t)\|_{H^{-1}(\R)} 
&\gtrsim \|q_{\alpha_c}(\cdot,0)\|_{H^{-1}(\R)} \gtrsim 
\alpha_c^{{\frac 12}}\|q(\cdot, 0)\|_{H^{-1}(\R)}.
\end{align*}
Recalling that $\alpha_c=R^{-2}$, we obtain
\[
R^{-2}\|q(\cdot,0)\|_{H^{-1}(\R)}\lesssim  \|q(\cdot,t)\|_{H^{-1}(\R)}\lesssim R^2\|q(\cdot,0)\|_{H^{-1}(\R)}\,
\]
for all $t \in \R$. } Finally, having proved \eqref{m1} for $q_0\in \Sch(\R)$, it is enough to use Theorem \ref{yap} to extend \eqref{m1} to $q_0\in L^2(\R)$.\qed

\smallskip

\noindent Our next goal is to prove the estimate
\begin{equation}\label{eq191}
C_1(1+\|q_0\|_{L^2(\R)})^{2s} \|q_0\|_{H^{s}(\R)}\le \|q(\cdot, t)\|_{H^{s}(\R)} \le C_2(1+\|q_0\|_{L^2(\R)})^{-2s} \|q_0\|_{H^{s}(\R)},
\end{equation} 
where $\,t \in \R\,, \,s\in (-1,0]$, and $C_1$ and $C_2$ are positive absolute constants.  For $s=0$, this bound is trivial. To cover $s\in (-1,0)$, we will need some auxiliary results first. One of the basic properties of NLS which we discussed in the Introduction has to do with modulation:  if $q(\xi,t)$ solves \eqref{nls1}, then $\widetilde q_v(\xi,t)=e^{iv\xi-iv^2t}q(\xi-2vt,t)$ solves \eqref{nls1} for every $v\in \R$.

\begin{Lem}\label{l52}
Let $q_0 \in L^2(\R)$, $t \in \R$. Then,
$$
\|\tilde q_v(\cdot, t)\|_{H^{-1}(\R)}^{2} =  \int_{\R}\frac{|(\F q)(\eta, t)|^2}{1+(\eta+v)^2}\,d\eta.
$$
\end{Lem}
\beginpf It is clear that  $\|e^{-iv^2t} f\|_{H^{-1}(\R)} = \|f\|_{H^{-1}(\R)}$ for every $f \in H^{-1}(\R)$ and $t \in \R$, because $e^{-iv^2t}$ is a unimodular constant. We have $\F (e^{iv\xi}q(\xi-2vt, t))(\eta) = (\F q(\xi, t))(\eta-v)e^{-2ivt(\eta-v)}$, $\eta \in \R$. Since $|e^{-2ivt(\eta-v)}| = 1$, it only remains to change the variable of integration in
$$
\|\widetilde q_v\|_{H^{-1}(\R)} = \int_{\R}\frac{|(\F q(\xi, t))(\eta-v)|^2}{1+\eta^2}\,d\eta
$$
to get the statement of the lemma.
\qed

\medskip

The next result is a standard property of convolutions.
\begin{Lem}\label{l53}
Let $\gamma \in (-\frac 12,1]$ and set $a_{k} = \frac{1}{(1+k^2)^{\gamma}}$ for  $k \in \Z$. We have
$$
 \sum_{k \in \Z} \frac{a_k}{1+(\eta-k)^2} \sim  C_\gamma \frac{1}{(1+\eta^2)^{\gamma}}, \qquad \eta \in \R.
$$
\end{Lem}
\beginpf After comparing  the sum to an integral, it is enough to show that
\[
\int_\R \frac{du}{(1+u^2)^\gamma(1+(\eta-u)^2)} \sim  C_\gamma\frac{1}{(1+\eta^2)^{\gamma}}\,.
\]
The function on the left-hand side is even and continuous in $\eta$ and $\gamma$, so we can assume that $\eta>1$. Then, 
\[
\int_{|\eta-u|<0.5\eta} \frac{du}{(1+u^2)^\gamma(1+(\eta-u)^2)} \sim  \frac{1}{(1+\eta^2)^{\gamma}},\quad
\int_{|\eta-u|>0.5\eta} \frac{du}{(1+u^2)^\gamma(1+(\eta-u)^2)} \lesssim \cal{I}_1+\cal{I}_2\,,
\]
where
\[
\cal{I}_1=\int_{u<-\eta/2,\; {u>3\eta/2}} \frac{du}{(1+u^2)^\gamma(1+(\eta-u)^2)}\lesssim \int_{-\infty}^{-\eta/2}\frac{du}{u^{2+2\gamma}}+\int_{3\eta/2}\frac{du}{u^{2+2\gamma}}   \le C_\gamma \eta^{-1-2\gamma}\,,
\]
\[
\cal{I}_2=\int_{|u|<\eta/2} \frac{du}{(1+u^2)^\gamma(1+(\eta-u)^2)}\lesssim  \eta^{-2}\int_{|u|<\eta/2} \frac{du}{(1+u^2)^\gamma}\lesssim \eta^{-2\gamma}\,.
\]
Combining these bounds proves the lemma.  \qed\medskip

\noindent{\bf Proof of Theorem \ref{t5}, the case $s \in (-1, 0)$.} We can again assume that $q_0\in \Sch(\R)$. Recall the estimate \eqref{eq19} for $s = -1$: 
\begin{equation}\label{lol1}
C_1(1+\|q_0\|_{L^2(\R)})^{-2} \|q_0\|_{H^{-1}(\R)}\le \|q(\cdot, t)\|_{H^{-1}(\R)} \le C_2(1+\|q_0\|_{L^2(\R)})^{2} \|q_0\|_{H^{-1}(\R)}\,.
\end{equation} 
According to Lemma \ref{l52}, we have
\begin{equation}\label{eq26}
\|\tilde q_v(\cdot, t)\|_{H^{-1}(\R)}^{2} =  \int_{\R}\frac{|(\F q(\cdot, t))(\eta)|^2}{1+(v+\eta)^2}\,d\eta
\end{equation}
for $\tilde q_{v}(\xi, t) = e^{iv\xi -iv^2t}q(\xi-2vt, t)$. Let $a_{k}$, $k \in \Z$, be the coefficients from Lemma \ref{l53} with $\gamma=-s$. Then, \eqref{eq26} and Lemma \ref{l53} imply
\begin{equation}\label{eq27}
 \sum_{k \in \Z}a_k
\|\tilde q_{k}(\cdot, t)\|_{H^{-1}(\R)}^{2}
\sim C_s
\|q(\cdot, t)\|_{H^{s}(\R)}^{2}\,.
\end{equation}
{In particular, taking $t=0$ gives
\begin{equation}\label{eq27bis}
\sum_{k \in \Z}a_k
\|\tilde q_{k}(\cdot, 0)\|_{H^{-1}(\R)}^{2}
\sim C_s
\|q_0\|_{H^{s}(\R)}^{2}\,.
\end{equation}}
We now apply \eqref{lol1} to $\tilde q_{k}$ {and use \eqref{eq27} and \eqref{eq27bis}} to get
\begin{equation}\label{eq192}
C_1(s)(1+\|q_0\|_{L^2(\R)})^{-2} \|q_0\|_{H^{s}(\R)}\le \|q(\cdot, t)\|_{H^{s}(\R)} \le C_2(s)(1+\|q_0\|_{L^2(\R)})^{2} \|q_0\|_{H^{s}(\R)}\,.
\end{equation} 
If $R=\|q_0\|_{L^2(\R)}\le 1$, we have the statement of our theorem. If $R=\|q_0\|_{L^2(\R)}> 1$, we  use dilation transformation like in the previous proof for $s=-1$. Consider $q_\alpha(\xi,t)=\alpha q(\alpha\xi, \alpha^2t)$ which solves the same equation and notice that
$
\|q_\alpha\|_{L^2(\R)}=\alpha^{\frac 12}R\,.
$
Let $\alpha=\alpha_c=R^{-2}<1$ making $\|q_{\alpha_c}\|_{L^2(\R)}=1$. Then, for the Sobolev norm, we have
\begin{equation}\notag
\|q_\alpha(\cdot,t)\|_{H^{s}(\R)}=\alpha^{\frac 12}\left(\int_\R \frac{1}{(1+\alpha^2\eta^2)^{|s|}}|(\cal{F} q)(\eta,{\alpha^2}t)|^2d\eta\right)^{\frac 12}\,.
\end{equation}
From \eqref{lop1}, 
\[
\frac{1}{(1+\eta^2)^{|s|}}\le \frac{1}{(1+\alpha_c^2\eta^2)^{|s|}}\le \frac{1}{\alpha_c^{2|s|}(1+\eta^2)^{|s|}}\,.
\]
Then, one has
\[
\alpha_c^{\frac 12}\|{q(\cdot, \alpha_c^2 t)}\|_{H^{s}(\R)}{\le} \|q_{\alpha_c}(\cdot,{t})\|_{H^{s}(\R)}{\le} \alpha_c^{\frac 12-|s|}\|{q(\cdot, \alpha_c^2 t)}\|_{H^{s}(\R)}\,.
\]
{In particular, taking $t = 0$ gives us
\[
\alpha_c^{\frac 12}\|q(\cdot, 0)\|_{H^{s}(\R)}
\le \|q_{\alpha_c}(\cdot,0)\|_{H^{s}(\R)}
\le \alpha_c^{\frac 12-|s|}\|q(\cdot, 0)\|_{H^{s}(\R)}\,.
\]
Now $\|q_{\alpha_c}(\cdot,0)\|_{L^2(\R)}=1$ and we can apply the previous bounds to get
\[
\|q_{\alpha_c}(\cdot,t)\|_{H^{s}(\R)}\sim  \|q_{\alpha_c}(\cdot,0)\|_{H^{s}(\R)}\,.
\]
Then,
\begin{align*}
\alpha_c^{\frac 12}\|q(\cdot, \alpha_c^2 t)\|_{H^{s}(\R)} 
&\lesssim \|q_{\alpha_c}(\cdot,0)\|_{H^{s}(\R)} \lesssim 
\alpha_c^{\frac 12-|s|}\|q(\cdot, 0)\|_{H^{s}(\R)},\\
\alpha_c^{\frac 12-|s|}\|q(\cdot, \alpha_c^2 t)\|_{H^{s}(\R)} 
&\gtrsim \|q_{\alpha_c}(\cdot,0)\|_{H^{s}(\R)} \gtrsim 
\alpha_c^{\frac 12}\|q(\cdot, 0)\|_{H^{s}(\R)}.
\end{align*}
Recalling that $\alpha_c=R^{-2}=\|q_0\|_{L^2(\R)}^{-2}$, we obtain
\[
\|q_0\|_{L^2(\R)}^{-2|s|}\|q(\cdot,0)\|_{H^{s}(\R)}\lesssim  \|q(\cdot,t)\|_{H^{s}(\R)}\lesssim \|q_0\|_{L^2(\R)}^{2|s|}\|q(\cdot,0)\|_{H^{s}(\R)}\,
\]
for all $t \in \R$.} \qed \smallskip

Our approach also provides the bounds for some positive Sobolev norms. The following proposition slightly improves \eqref{koch} when $s\in [0,\frac 12)$,  { $\|q_0\|_{H^s(\R)}$ is large, and $\|q_0\|_{L^2(\R)}$ is much smaller than $\|q_0\|_{H^s(\R)}$.}
\begin{Prop}\label{p5}
Let $q_0 \in \Sch(\R)$  and let $q = q(\xi, t)$ be the solution of \eqref{nls} corresponding to $q_0$. Then,  for each $s\in [0,\frac 12)$, we get
\begin{equation}\label{nes0}
 \|q(\cdot, t)\|_{H^{s}(\R)} \sim C_s \|q_0\|_{H^{s}(\R)}
\end{equation}
if $\|q_0\|_{L^2(\R)}\le 1$ and 
\begin{equation}\label{nes1}
 \|q(\cdot, t)\|_{H^{s}(\R)} \lesssim  C_s(\|q_0\|_{L^2(\R)}^{1+2s}+ \|q_0\|_{H^{s}(\R)})
\end{equation}
if $\|q_0\|_{L^2(\R)}>1$.
\end{Prop}

\beginpf In the case when $\|q\|_{L^2(\R)}\le 1$, the proof of proposition repeats the arguments given above to get \eqref{eq192} except that the constants in the inequalities depend on $s$ and can blow up when $s\to \frac 12$.
Suppose $\|q\|_{L^2(\R)}\ge 1$. Then, for the Sobolev norm, we have
\begin{equation}\notag
\|q_\alpha(\cdot,t)\|_{H^{s}(\R)}=\alpha^{\frac 12}\left(\int_\R {(1+\alpha^2\eta^2)^{s}}|(\cal{F} q)(\eta,{\alpha^2}t)|^2d\eta\right)^{\frac 12}\,.
\end{equation}
Take $\alpha=\alpha_c$ and write the following estimate for the integral above:
\begin{eqnarray*}
\int_\R {\left(1+\frac{\eta^2}{R^4}\right)^{s}}|(\cal{F} q)(\eta,{\alpha_c^2}t)|^2d\eta\sim \int_{-R^2}^{R^2} |(\cal{F} q)(\eta,{\alpha_c^2}t)|^2d\eta+R^{-4s}\int_{|\eta|>R^2} {\left(1+\eta^2\right)^{s}}|(\cal{F} q)(\eta,{\alpha_c^2}t)|^2d\eta\\
\lesssim R^2+R^{-4s}\int_{\R} {\left(1+\eta^2\right)^{s}}|(\cal{F} q)(\eta,{\alpha_c^2}t)|^2d\eta\, .\notag
\end{eqnarray*}
We use $\|q_{\alpha_c}(\cdot,t)\|_{L^2(\R)}=1$ and \eqref{eq192} to get $\|q_{\alpha_c}(\cdot,t)\|_{H^s(\R)}\sim C_s\|q_{\alpha_c}(\cdot,0)\|_{H^s(\R)}$.  The previous estimate for $t=0$ yields
$
\|q_{\alpha_c}(\cdot,0)\|_{H^s(\R)}\lesssim 1+R^{-1-2s}\|q(\cdot,0)\|_{H^s(\R)}.
$
Hence, 
$
\|q_{\alpha_c}(\cdot,t)\|_{H^s(\R)}\le C_s(1+R^{-1-2s}\|q(\cdot,0)\|_{H^s(\R)}).
$
We can write a lower bound
\[
\|q_{\alpha_c}(\cdot,t)\|_{H^{s}(\R)}^2=
R^{-2}\int_\R {\left(1+\frac{\eta^2}{R^4}\right)^{s}}|(\cal{F} q)(\eta,{\alpha_c^2}t)|^2d\eta\gtrsim R^{-2-4s}\int_{|\eta|>R^2} {\left(1+\eta^2\right)^{s}}|(\cal{F} q)(\eta,{\alpha_c^2}t)|^2d\eta
\]
so
\[
 \int_{|\eta|>R^2} {\left(1+\eta^2\right)^{s}}|(\cal{F} q)(\eta,{\alpha_c^2}t)|^2d\eta\le C_s(R^{2+4s}+\|q(\cdot,0)\|^2_{H^s(\R)}).
\]
Writing the integral as a sum of two:
\[
\int_{\R} {\left(1+\eta^2\right)^{s}}|(\cal{F} q)(\eta,{\alpha_c^2}t)|^2d\eta=\int_{|\eta|>R^2} {\left(1+\eta^2\right)^{s}}|(\cal{F} q)(\eta,{\alpha_c^2}t)|^2d\eta+\int_{|\eta|<R^2} {\left(1+\eta^2\right)^{s}}|(\cal{F} q)(\eta,{\alpha_c^2}t)|^2d\eta\,
\]
and estimating each of them, we get a bound which holds for all $t$:
\[
\int_{\R} {\left(1+\eta^2\right)^{s}}|(\cal{F} q)(\eta,{\alpha_c^2}t)|^2d\eta\le C_s(R^{2+4s}+\|q(\cdot,0)\|^2_{H^s(\R)})\,.
\]
That is the required upper bound \eqref{nes1}. \qed

\bigskip

\section{Oscillation and Sobolev space $H^{-1}(\R)$.}

In this part of the paper, our goal is to prove the Theorem \ref{P2c1b}. Let us recall its statement.
\begin{Thm}\label{tt}
Suppose that $q\in L^2(\R)$ and let $N_Q$ satisfy $JN'_Q + QN_Q = 0, N_Q(0)=I$, where  $Q = \left(\begin{smallmatrix}-\Im q & -\Re q \\ -\Re q & \Im q\end{smallmatrix}\right).$
Then, 
\begin{eqnarray}
e^{-C_1R}\|q\|_{H^{-1}(\R)}^2\lesssim \widetilde \K_{Q}\lesssim  e^{C_2R}\|q\|^2_{H^{-1}(\R)}\,,
\end{eqnarray}
where $R=\|q\|_{L^2(\R)}$ and $C_1$, $C_2$ are two positive absolute constants.
\end{Thm}

 {Theorem \ref{tt} is of independent interest in the spectral theory of Dirac operators.  For example, Lemma \ref{l88} shows  that $\|q\|_{L^2(\R)}$ and $\|q\|_{H^{-1}(\R)}$ control the size of $\K_Q$.}

\medskip 

The strategy of the proof is the following. In the next subsection, we show that $H^{-1}(\R)$-norm of any function can be characterized through BMO-like  condition for its ``antiderivative''. In Subsection~\ref{P2sec2}, we consider solution to Cauchy problem $JN'+QN=0, N(0)=I$ on the interval $[0,1]$ where zero-trace symmetric $Q$  and study the quantity $\det \int_0^1 N^*Ndx$, which represents a single term in the sum for $\widetilde \K_Q$. The results in Subsection~\ref{P2sec3} show that small value of $\widetilde \K_Q$ guarantees that the ``local'' $H^{-1}$ norm of $Q$ is also small. This rough estimate is used in the proof of Theorem \ref{P2c1b} which is contained in~Subsection~\ref{P2sec4}.
\medskip

\subsection{One property of Sobolev space $H^{-1}(\R)$}\label{P2sec1} Observe that a function $f \in L^2(\R)$ belongs to the Sobolev space $H^{-1}(\R)$ if and only if
\begin{equation}\label{eq1}
 \int_{\R}\left|\int_{\R}f(y)\chi_{\R_+}(x-y)e^{-(x-y)}\,dy\right|^2\,dx < \infty.
\end{equation}
Moreover, the last integral is equal to $\|f\|_{H^{-1}(\R)}^{2}.$
Indeed, recall that $\F f$ stands for the Fourier transform of $f$:
$$
(\F f)(\eta) = \frac{1}{\sqrt{2\pi}}\int_{\R}f(x)e^{-i\eta x}\,dx\,.
$$
Then, from Plancherel's identity and formula
$$
\frac{1}{\sqrt{2\pi}}\int_{\R_+}e^{-x}e^{-ix\eta}dx = \frac{1}{\sqrt{2\pi}}\frac{1}{1+i\eta},
$$
we obtain
$$
\|f\|_{H^{-1}(\R)}^{2} = 2\pi\|(\F f ) \F(\chi_{\R_+}e^{-x})\|^2_{L^2(\R)} = \int_{\R}\frac{|(\F f)(\eta)|^2}{1+\eta^2}\,d\eta
$$
by properties of convolutions.
We will need the following  proposition. 

\begin{Prop}\label{P2p1lol}
Suppose that $f \in L^1_{\rm  loc}(\R) \cap H^{-1}(\R)$. Let $g$ be an absolutely continuous function on $\R$ such that $g' = f$ almost everywhere on $\R$. Then,
\begin{equation}\label{P2eq4}
c_1\|f\|_{H^{-1}(\R)}^{2} \le \sum_{k \in \Z} \int_{I_{k}}|g - \langle g \rangle_{I_{k}}|^2\,dx \le c_2\|f\|_{H^{-1}(\R)}^{2},
\end{equation}
where $I_{k} = [k, k+2]$, $\langle g \rangle_{I} = \frac{1}{|I|}\int_{I} g(x)\,dx$, and the positive constants $c_{1}$ and $c_2$ are universal.
\end{Prop}
\beginpf Take a function $f \in L^1_{\rm  loc}(\R) \cap H^{-1}(\R)$, and let $g$ be an absolutely continuous function on $\R$ such that $g' = f$ almost everywhere on $\R$.  Assume first that $f$ has a compact support. The integral under the sum does not change if we add a constant to $g$, so we can suppose without loss of generality that
$$
g(x) = \int_{-\infty}^{x}f(s)\,ds, \qquad x \in \R.
$$ 
{\bf Upper bound.} 
Given $f$, define {$o_f$} by
\[
o_f(x)=e^{-x}\int_{-\infty}^x f(y)e^ydy
\]
and recall (see \eqref{eq1}) that 
\begin{equation}\label{sp1}
\|f\|_{H^{-1}(\R)}=\|o_f\|_{L^2(\R)}.
\end{equation}
Moreover, 
\begin{equation}
o'_f+o_f=f\,. \label{P1sa1}
\end{equation}
For each interval $I_k$, we use \eqref{P1sa1} for the corresponding term in the sum \eqref{P2eq4}:
\begin{eqnarray*}
\int_{I_k} \left|\int_k^x fdx_1-\frac 12\int_k^{k+2} \left(\int_k^{x_1}f(x_2)dx_2\right)dx_1\right|^2dx\\
=\int_{I_k} \left|\int_k^x o(x_1)dx_1+o(x)-\frac 12\int_k^{k+2} \left(o(x_1)+\int_k^{x_1}o(x_2)dx_2\right)dx_1\right|^2dx\\
\lesssim \int_{I_k}|o|^2dx
\end{eqnarray*}
after the Cauchy-Schwarz inequality is applied. Summing these estimates in $k\in \Z$ and using \eqref{sp1}, we get the upper bound in \eqref{P2eq4} for compactly supported $f$. \smallskip

\noindent {\bf Lower bound.} Integration by parts gives
\begin{align*}
\int_{-\infty}^{x}f(y)e^{-(x-y)}\,dy 
&= \int_{-\infty}^{x}f(s)\,ds - \int_{-\infty}^{x}\left(\int_{-\infty}^{y}f(s)\,ds\right) \, e^{-(x-y)}\,dy\\
&= g(x) - \int_{-\infty}^{x}g(y)\,e^{-(x-y)}\,dy\\ 
&= 
\int_{-\infty}^{x}(g(x)-g(y))\,e^{-(x-y)}\,dy.
\end{align*}
Therefore, 
\begin{align*}
\int_{\R}\left|\int_{-\infty}^{x}f(y)e^{-(x-y)}\,dy\right|^2\,dx 
&\le \sum_{k \in \Z} \int_{k}^{k+2}\left(\int_{-\infty}^{x}|g(x)-g(y)|^2e^{-(x-y)}\,dy\right)\,dx\\
&\lesssim \sum_{k \in \Z}\sum_{j \le k}e^{-(k-j)} \int_{k}^{k+2}\int_{j}^{j+2}|g(x) - g(y)|^2\,dx\,dy\,.   
\end{align*}  
Using the inequality $(x+y+z)^2 \le 3(x^2+y^2+z^2)$, we continue the estimate:
$$
... \lesssim \sum_{k \in \Z}\sum_{j \le k}e^{-(k-j)}\left(\int_{I_k}|g - \langle g \rangle_{I_k}|^2dx + |\langle g \rangle_{I_j} - \langle g \rangle_{I_k}|^2 + \int_{I_{j}}|g - \langle g \rangle_{I_j}|^2dx\right). 
$$
Since
$$
\sum_{k \in \Z}\sum_{j \le k}e^{-(k-j)}\left(\int_{I_k}|g - \langle g \rangle_{I_k}|^2dx + \int_{I_{j}}|g - \langle g \rangle_{I_j}|^2dx\right) \lesssim \sum_{k\in \Z} \int_{I_k}|g - \langle g \rangle_{I_k}|^2dx\,,
$$
we are left with estimating
$$
\sum_{k \in \Z}\sum_{j \le k}e^{-(k-j)}|\langle g \rangle_{I_j} - \langle g \rangle_{I_k}|^2.
$$
Applying the Cauchy-Schwarz inequality for the telescoping sum
\[
\langle g\rangle_{I_k}-\langle g\rangle_{I_j}=\sum_{s=j+1}^k\Bigl(\langle g\rangle_{I_s}-\langle g\rangle_{I_{s-1}}\Bigr)\,,
\]
we get
$$
|\langle g \rangle_{I_j} - \langle g \rangle_{I_k}|^2 \le (k-j) \sum_{j\le s \le k-1}|\langle g \rangle_{I_s} - \langle g \rangle_{I_{s+1}}|^2. 
$$ 
Then,  
$$
\sum_{k \in \Z}\sum_{j \le k}e^{-(k-j)}(k-j) \sum_{j\le s \le k-1}|\langle g \rangle_{I_s} - \langle g \rangle_{I_{s+1}}|^2 = \sum_{s \in \Z} |\langle g \rangle_{I_s} - \langle g \rangle_{I_{s+1}}|^2 \sum_{k, j: \; j \le s \le k-1} (k-j)e^{-(k-j)}.
$$
We have
\begin{align*}
\sum_{k, j: \; j \le s \le k-1} (k-j)e^{-(k-j)}
&=\sum_{j \le s}\sum_{m\ge 1} (s+m-j)e^{-(s+m-j)} \\
&= \sum_{j \le 0}\sum_{m\ge 1} (m-j)e^{-(m-j)} = \sum_{j \ge 0}\sum_{m\ge 1} (m+j)e^{-(m+j)}.  
\end{align*}
The last sum is finite and does not depend on index $s$. Now, the estimate 
$$
|\langle g \rangle_{I_s} - \langle g \rangle_{I_{s+1}}|^2 = \int_{I_{s} \cap I_{s+1}}|\langle g \rangle_{I_s} - g + g - \langle g \rangle_{I_{s+1}}|^2dx \le
2\int_{I_s}|g - \langle g \rangle_{I_s}|^2dx + 2\int_{I_{s+1}}|g - \langle g \rangle_{I_{s+1}}|^2dx
$$
proves that 
$$
\sum_{k \in \Z}\sum_{j \le k}e^{-(k-j)}|\langle g \rangle_{I_j} - \langle g \rangle_{I_k}|^2\lesssim \sum_{s\in \Z} \int_{I_s}|g - \langle g \rangle_{I_s}|^2dx\,.
$$
Hence, the lower bound in \eqref{P2eq4} holds for compactly supported $f$. 

\smallskip

Now, take any $f\in  L^1_{\rm  loc}(\R) \cap H^{-1}(\R)$. The definition \eqref{n1} of $H^{-1}(\R)$ implies that $\F f$ can be written as $(1+i\eta)(\F o)$ for some function $o\in L^2(\R)$. Moreover, this map $f\mapsto o$ is a bijection between $H^{-1}(\R)$ and $L^2(\R)$ and $\|f\|_{H^{-1}(\R)}=\|o\|_{L^2(\R)}$. Taking the inverse Fourier transform of identity $\F f=(1+i\eta)(\F o)$, one gets a formula $f=o+o'$ where $o'$ is understood as a derivative in $\Sch'(\R)$. Since $f\in L^1_{\rm loc}(\R)$ and $o\in L^2(\R)$, we have $o'\in L^1_{\rm loc}(\R)$ and, therefore, $o$ is absolutely continuous on $\R$ with the derivative equal to $f-o$. Now, take $o_n(x)=o(x)\mu_n(x)$ and define the corresponding $f_n=o_n+o_n'$. Here, $\mu_n(x)$ is even and 
\[
\mu_n(x)=\left\{  
\begin{array}{ll}
1, &0\le x<n,\\
n+1-x,& x\in [n,n+1),\\
0, &x\ge n+1\,.
\end{array}\right.
\]
Then, $\{o_n\}\to o$ in $L^2(\R)$ and so $\{f_n\}\to f$ in $H^{-1}(\R)$  {because the mapping $f \mapsto o$ is unitary from $H^{-1}(\R)$ onto $L^2(\R)$}.  Also, each $f_n$ is compactly supported and $\{f_n\}$ converges to $f$ uniformly on every finite interval. Define $g_n=\int_0^x f_nds, g=\int_0^x fds$, and write  \eqref{P2eq4} for $f_n$.  The estimate on the right gives
\[
 \sum_{|k|\le N} \int_{I_{k}}|g_n - \langle g_n \rangle_{I_{k}}|^2\,dx \le c_2\|f_n\|_{H^{-1}(\R)}^{2}
\]
for each $N\in \mathbb{N}$. Sending $n\to \infty$, the bound
\[
 \sum_{|k|\le N} \int_{I_{k}}|g - \langle g\rangle_{I_{k}}|^2\,dx \le c_2\|f\|_{H^{-1}(\R)}^{2}
\]
appears.
Taking $N\to\infty$, one has the right estimate in \eqref{P2eq4}. In particular, it shows that the sum in \eqref{P2eq4} converges. By construction, 
\[
 \sum_{k\in \Z} \int_{I_{k}}|g_n - \langle g_n \rangle_{I_{k}}|^2\,dx=
 \sum_{-n\le k\le n-2} \int_{I_{k}}|g - \langle g \rangle_{I_{k}}|^2\,dx+\epsilon_n\,,
\]
where $\epsilon_n$ is a sum of integrals over $I_{-n-2},I_{-n-1},I_{n-1},I_n$. Since $o\in L^2(\R)$, 
\[
\lim_{n\to\infty}\int_{I_{k}}|g_n - \langle g_n \rangle_{I_{k}}|^2\,dx=0, \quad k\in \{-n-2,-n-1,n-1,n\}\,.
\]
Hence, $\lim_{n\to\infty}\epsilon_n=0$ and, taking $n\to\infty$ in inequality
\[
c_1\|f_n\|_{H^{-1}(\R)}^{2} \le \sum_{k \in \Z} \int_{I_{k}}|g_n - \langle g_n \rangle_{I_{k}}|^2\,dx\,,
\]
one gets the left bound in \eqref{P2eq4}. Since all antiderivatives are different by a constant and the integral in \eqref{P2eq4} does not change if we add a constant to $g$, the proof is finished.
\qed

\medskip

\medskip

\subsection{Auxiliary perturbative results for a single interval}\label{P2sec2}

 Notice that for any real symmetric $2\times 2$ matrix $Q$ with zero trace, we have that $V=JQ$ is also real, symmetric and has zero trace. The converse statement is true as well. Hence, the equation $JN_Q' + QN_Q = 0$  in Theorem \ref{tt}, which is equivalent to $N_Q'=JQN_Q$, can be written as $N_Q'=VN_Q$ with $V$ having the same properties as $Q$. Let $U_+(x,y)$ denote the solution to \[\frac{d}{dx}U_+(x,y)=V(x)U_+(x,y),\quad U_+(y,y)=I\] and $U_-(x,y)$ denote the solution to \[\frac{d}{dx}U_-(x,y)=-V(x)U_-(x,y), \quad U_-(y,y)=I\,.\] 
\begin{Lem}Suppose $N'=VN,  N(0)=I$, where $V$ is real-valued, $V\in L^1[0,1]$, $V=V^*$,  and ${\rm tr}\,V=0$. Then, for $\Hh=N^*N$, we have
\begin{eqnarray}\label{uo}
\det \int_0^1 \Hh(\xi)\,d\xi=\frac 12\int_0^1\int_0^1 {\rm tr}\, \Bigl(U^*_+(x,y)U_+(x,y)\Bigr)\,dxdy=
\frac 12\int_0^1\int_0^1 \|U_+(x,y)\|^2_{\rm HS}\,dxdy,\\
\det \int_0^1 \Hh(\xi)\,d\xi-1=
\frac 12\int_0^1\int_0^1 \left\|\Bigl(U_+(x,y)-U_-(x,y)\Bigr)e_1\right\|^2dxdy\,.\label{uu}
\end{eqnarray}
\end{Lem}
\beginpf Notice that $N, U_+,U_-\in {\rm SL}(2,\R)$ and that
every matrix $A\in {\rm SL}(2,\R)$ satisfies
\begin{equation}\label{py}
JA^*=A^{-1}J,\qquad AJ=J(A^*)^{-1}\,.
\end{equation}
Also, for any real $2\times 2$ matrix $B$, we have
\[
\det B=\langle JBe_1,Be_2\rangle=-\langle JBe_2,Be_1\rangle\,.
\]
Hence, 
\begin{eqnarray*}
\mathcal{I}:=\det \int_0^1 \Hh(\xi)d\xi=\int_0^1\int_0^1 \langle JN^*(x)N(x)e_1, N^*(y)N(y)e_2\rangle dxdy\\=-\int_0^1\int_0^1 \langle JN^*(x)N(x)e_2, N^*(y)N(y)e_1\rangle dxdy\,.
\end{eqnarray*}
For the second integrand, we have
\[
\langle JN^*(x)N(x)e_1, N^*(y)N(y)e_2\rangle=\langle N^*(y)N(y)JN^*(x)N(x)e_1, e_2\rangle\,.
\]
Then, identities \eqref{py} imply
\[
N^*(y)N(y)JN^*(x)N(x)=N^*(y)J(N^*(y))^{-1}N^*(x)N(x)=J(N(y))^{-1}(N^*(y))^{-1}N^*(x)N(x)
\]
and, since $Je_1=e_2$ and $J^*=-J$,
\[
\langle JN^*(x)N(x)e_1, N^*(y)N(y)e_2\rangle=\langle (N(y))^{-1}(N^*(y))^{-1}N^*(x)N(x)e_1,e_1\rangle\,.
\]
Similarly, 
$
\langle JN^*(x)N(x)e_2, N^*(y)N(y)e_1\rangle=-\langle (N(y))^{-1}(N^*(y))^{-1}N^*(x)N(x)e_2,e_2\rangle\,.
$
Hence,
\begin{eqnarray*}
\mathcal{I}=\frac 12 \int_0^1\int_0^1 \sum_{j=1}^2 \langle (N(y))^{-1}(N^*(y))^{-1}N^*(x)N(x)e_j,e_j\rangle dxdy=\hspace{3cm}\\
\frac 12 \int_0^1\int_0^1 {\rm tr} \Bigl((N(y))^{-1}(N^*(y))^{-1}N^*(x)N(x)\Bigr) dxdy=\frac 12 \int_0^1\int_0^1 {\rm tr} \Bigl((N^*(y))^{-1}N^*(x)N(x)(N(y))^{-1}\Bigr)  dxdy\,.
\end{eqnarray*}
Now, we use the formula $N(x)(N(y))^{-1}=U_+(x,y)$  to rewrite the last expression as
\[
\mathcal{I}=\frac 12 \int_0^1\int_0^1 {\rm tr} \Bigl(U^*_+(x,y)U_+(x,y)\Bigr)  dxdy\,.
\]
Finally,  \eqref{uu} follows from $U_+(x,y)\in {\rm SL}(2,\R)$ by direct inspection after one uses the identities
$JU_+(x,y)J=-U_-(x,y)$ and $ {\rm tr}(A^*A)-2=\|(A+JAJ)e_1\|^2$, which holds for every $A\in {\rm SL}(2,\R)$.
\qed\smallskip

\noindent{\bf Remark.} The integrand in \eqref{uo} is symmetric: $ {\rm tr} \Bigl(U^*_+(x,y)U_+(x,y)\Bigr)= {\rm tr} \Bigl(U^*_+(y,x)U_+(y,x)\Bigr)$ because $U_+(x,y)=U^{-1}_+(y,x)$ and $U_+(x,y)\in {\rm SL}(2,\R)$. Notice also, that 
\[
 {\rm tr} \Bigl(U^*_+(x,y)U_+(x,y)\Bigr)=\lambda_{x,y}^2+\lambda_{x,y}^{-2}\ge 2\,,
\]
where $\lambda_{x,y}$ is an eigenvalue of $U^*_+(x,y)U_+(x,y)$ which explains why the left-hand side in \eqref{uu} is nonnegative.

\begin{Lem}\label{P2le0}
 Suppose real-valued  matrix-function $V=
\left(\begin{smallmatrix}
v_1 & v_2\\
v_2 & -v_1
\end{smallmatrix}
\right)
$ 
is  defined on $[0,1]$ and
satisfies $\|V\|_{L^1[0,1]}<\infty$. Consider $\Hh=N^*N$, where $N$: $N'=VN, N(0)=I$.
Then, 
\begin{equation}\label{P2oi1}
\det \int_0^1 \Hh\,dx-1\lesssim  \|V\|^2_{L^1[0,1]}\exp 
(C\|V\|_{L^1[0,1]})\,.
\end{equation}
\end{Lem}
\beginpf
The integral equation for $N$ is
\begin{equation}\label{P2op1}
N=I+\int_0^x VNds\,.
\end{equation}
By Gronwall's inequality, 
\begin{equation}\label{lok}
\|N(x)\|\le \exp\left(\int_0^x\|V(s)\|ds\right)\le \exp (\|V\|_{L^1[0,1]}).
\end{equation}
 Iteration of \eqref{P2op1} gives
\[
N=I+\int_0^x Vdx_1+\int_0^xV(x_1)\left(\int_0^{x_1}V(x_2)N(x_2)dx_2\right)dx_1\,.
\]
Then, 
\[
\int_0^1 N^*Ndx=I+2\int_0^1 \left(\int_0^xV(x_1)dx_1\right)dx+R, \;\; \|R\|\lesssim \|V\|^2_{L^1[0,1]}\exp (C\|V\|_{L^1[0,1]})\,.
\]
Since $\text{tr} \,V=0$, the identity
$
\det (I+A)=1+\text{tr} A+\det A\,,
$
which holds for all $2\times 2$ matrices $A$,
gives
\[
\det \int_0^1 \Hh dx-1\lesssim \|V\|^2_{L^1[0,1]}\exp (C\|V\|_{L^1[0,1]})\,.
\]
\qed

\begin{Lem}\label{P2le1} Suppose real-valued symmetric matrix-functions $V$ and $O$ are defined on $[0,1]$ and satisfy
\begin{eqnarray}
V=\left(
\begin{matrix}
v_1 & v_2\\
v_2 & -v_1
\end{matrix}
\right)=
O+O', O=O^*=\left(
\begin{matrix}
o_1 & o_2\\
o_2 & -o_1
\end{matrix}
\right)\,,\\
\delta:=\|O\|_{L^2[0,1]}<\infty\,,   \label{P2p1}\\
d:=\|O'\|_{L^2[0,1]}<\infty\,.\label{P2p4}
\end{eqnarray}
Consider $\Hh=N^*N$ where $N'=VN, N(0)=I$.
Then, we have
\begin{equation}\label{resid}
\det \int_0^1 \Hh \, dx-1 = 4\sum_{j=1}^2 \int_0^1 |g_j-\langle g_j\rangle|^2\,dx+r, \, |r|\lesssim \delta^{2.5}\exp(C(d+\delta))\,,
\end{equation}
where 
\begin{equation}\label{okl1}
g_j:=\int_0^x v_j\,dx
\end{equation}
and $C$ is an absolute positive constant. An analogous result holds if $O$ and $V$ are related by $V=O-O'$.
\end{Lem}
\beginpf 
We will use the formula \eqref{uu} for our analysis. Fix $y\in [0,1]$ and take $U_+(x,y)$ and $U_-(x,y)$ which solve $\frac{d}{dx}U_+(x,y)=V(x)U_+(x,y), U_+(y,y)=I$ and $\frac{d}{dx}U_-(x,y)=-V(x)U_-(x,y), U_-(y,y)=I$. Iterating the corresponding integral equations, one gets
\begin{eqnarray*}\nonumber
U_+(x,y)=I+\int_y^x Vdx_1+\int_y^x V\int_y^{x_1}Vdx_2dx_1+\int_y^x V\int_y^{x_1}V\int_y^{x_2}Vdx_3dx_2dx_1+\\
\int_y^x V\int_y^{x_1}V\int_y^{x_2}V\int_y^{x_3}Vdx_4dx_3dx_2dx_1+\int_y^x V\int_y^{x_1}V\int_y^{x_2}V\int_y^{x_3}Vf_+dx_4dx_3dx_2dx_1,\\ f_+(x_4)=\int_y^{x_4} V(s)U_+(s,y)ds\,.\nonumber
\end{eqnarray*}
\begin{eqnarray*}
U_-(x,y)=I-\int_y^x Vdx_1+\int_y^x V\int_y^{x_1}Vdx_2dx_1-\int_y^x V\int_y^{x_1}V\int_y^{x_2}Vdx_3dx_2dx_1+\\
\int_y^x V\int_y^{x_1}V\int_y^{x_2}V\int_y^{x_3}Vdx_4dx_3dx_2dx_1-\int_y^x V\int_y^{x_1}V\int_y^{x_2}V\int_y^{x_3}Vf_-dx_4dx_3dx_2dx_1,\\  f_-(x_4)=\int_y^{x_4} V(s)U_-(s,y)ds\,.\nonumber
\end{eqnarray*}
Taking $U_+(x,y)-U_-(x,y)$ as in \eqref{uu}  leaves us with 
\begin{eqnarray}\label{pip}
\frac{U_+(x,y)-U_-(x,y)}{2}=\int_y^x Vdx_1+\cal{I}_1+\cal{I}_2\,,\\
\cal{I}_1=\int_y^x V\int_y^{x_1}V\int_y^{x_2}Vdx_3dx_2dx_1,\\
\cal{I}_2=
\int_y^x V\int_y^{x_1}V\int_y^{x_2}V\int_y^{x_3}V(f_++ f_-)dx_4dx_3dx_2dx_1\,.
\end{eqnarray}
Recall that $V=O+O'$ where $O$ satisfies \eqref{P2p1} and \eqref{P2p4}. These assumptions are to be used in the following proposition. On $\R_+^2$, we define the partial order 
\[
\left[
\begin{array}{c}
x_1\\x_2
\end{array}
\right]\le \left[
\begin{array}{c}
y_1\\y_2
\end{array}
\right]
\]
by requiring that  $x_1\le y_1$ and $x_2\le y_2$.
\begin{Prop}
Suppose a matrix-function $O$ is defined on $[0,1]$ and denote 
\begin{eqnarray}
\delta=\|O\|_{L^2[0,1]}, \,
d=\|O'\|_{L^2[0,1]}.
\end{eqnarray}
Let an operator $G_{(y)}$ be given by: $F\mapsto (G_{(y)}F)(x)=\int_y^x (O+O')Fds$ where  $y\in [0,1]$ and a matrix-function $F$, defined on $[0,1]$, satisfies $\|F\|_{L^\infty[0,1]}<\infty$ and $\|F'\|_{L^2[0,1]}<\infty$. Then,
\begin{equation}\label{lop}
\left[
\begin{array}{l}
 \|G_{(y)}F\|_{L^\infty[0,1]}\\
 \|(G_{(y)}F)'\|_{L^2[0,1]}
\end{array}
\right]\le C \cal{M}
 \left[
\begin{array}{c}
\|F\|_{L^\infty[0,1]}
\\
\|F'\|_{L^2[0,1]}
\end{array}
\right]\,,\quad \cal{M}=\left(
\begin{array}{ll}
\delta+\sqrt{\delta d} & \delta \\
\delta+d &0
\end{array}
\right)\,,
\end{equation}
where $C$ is an absolute positive constant, the norms and derivatives are computed with respect to $x$.
\end{Prop}
\beginpf
 Let $b=\|F\|_{L^\infty[0,1]}, c=\|F'\|_{L^2[0,1]}$.
Write 
\begin{equation}\label{lak}
O^*(x)O(x)-O^*(y)O(y)=\int_y^x ((O^*)'O+O^*O')ds\,.
\end{equation}
Then, 
\[
\|O(x)\|^2=\max_{\|\xi\|_{\C^2}\le 1}\|O(x)\xi\|^2_{\C^2}=\max_{\|\xi\|_{\C^2}\le 1}\langle O^*(x)O(x)\xi,\xi\rangle\stackrel{\eqref{lak}}{\le} \|O(y)\|^2+2\int_0^1 \|O'(s)\|\cdot \|O(s)\|ds\,.
\]
Applying Cauchy-Schwarz inequality to the integral, integrating in $y$ from $0$ to $1$ and maximizing in $x$ gives
\begin{equation}
\|O\|_{L^\infty[0,1]}\lesssim \delta+(d\delta)^{\frac 12}\,. \label{P2p2}
\end{equation}
Then, 
\[
(G_{(y)}F)(x)=\int_y^x OFds+O(x)F(x)-O(y)F(y)-\int_y^xOF'ds
\]
and the estimate for the first coordinate in \eqref{lop} follows from Cauchy-Schwarz inequality and \eqref{P2p2}. Since 
$
(G_{(y)}F)'=(O+O')F\,,
$
we get $\|(G_{(y)}F)'\|_{L^2[0,1]}\le (\|O\|_{L^2[0,1]}+\|O'\|_{L^2[0,1]})\|F\|_{L^\infty[0,1]}$ and the bound for the second coordinate in \eqref{lop} is obtained.
\qed
\smallskip

\noindent{\bf Continuation of the proof of Lemma \ref{P2le1}.} We apply the proposition to $\cal{I}_1$ three times with the initial  choice of $F$: $F=I$. That gives rise to taking the third power of matrix $\mathcal{M}$: $\mathcal{M}^3$, applying it to $(1,0)^t$, and looking at the first coordinate. As the result, one has  $\|\cal I_1\|_{L^\infty[0,1]}\lesssim \delta^{\frac 32}(\delta+d)^{\frac 32}$. Therefore, 
\begin{equation}\label{eas}
\|\cal{I}_1e_1\|_{L^\infty([0,1]^2)} \lesssim \delta^{\frac 32}\exp(\delta+d).
\end{equation} Similarly, we consider $\cal{I}_2$ and use the previous proposition four times making the first choice of $F$ as $F=f_++f_-$. Applying the bound \eqref{lok} to $U_+$ and $U_-$, we get $ \|f_++f_-\|_{L^\infty[0,1]}\lesssim (\delta+d)\exp(\delta+d), \|f'_++f_-'\|_{L^2[0,1]}\lesssim (\delta+d)\exp(\delta+d)$. This time, we compute the fourth power of matrix $\mathcal{M}: \mathcal{M}^4$, apply it to vector $(\delta+d)\exp(\delta+d)(1,1)^t$, and look at the first coordinate. 
In the end, one has 
\begin{equation}\label{eas0}
\|\cal{I}_2e_1\|_{L^\infty([0,1]^2)}\lesssim \delta^{2}\exp(C(d+\delta))\,.
\end{equation} The first term in \eqref{pip} can be written as
\[
\int_y^x Vds=\int_y^x Ods+O(x)-O(y)
\]
and 
\begin{equation}\label{eas1}
\left\|\int_y^x Ods+O(x)-O(y)\right\|_{L^2([0,1]^2)}\lesssim  \delta\,.
\end{equation}
For any three vectors $v_1,v_2$ and $v_3$ in $\R^2$, we have an estimate
 {
\[
|\|v_1+v_2+v_3\|-\|v_1\|| 
\le \|v_2+v_3\| \le \|v_2\|+\|v_3\|\,,
\]
which follows from the triangle inequality.
Multiplying with 
\[
\|v_1+v_2+v_3\|+\|v_1\|
\le 2\|v_1\|+\|v_2\|+\|v_3\|\,,
\]
we get
\[
|\|v_1+v_2+v_3\|^2-\|v_1\|^2|
\le 2\|v_1\|(\|v_2\|+\|v_3\|)+(\|v_2\|+\|v_3\|)^2\,.
\]
}
Applying it to \eqref{pip} gives
\begin{eqnarray*}
\left|\frac 14 \|(U_+(x,y)-U_-(x,y))e_1\|^2-\left\|\left(\int_y^x Vds\right)e_1\right\|^2\right|   \hspace{4cm}\\
\lesssim \left\|\left(\int_y^x Vds\right)e_1\right\|\cdot (\|\cal{I}_1e_1\|+\|\cal{I}_2e_1\|)+\|\cal{I}_1e_1\|^2+\|\cal{I}_2e_1\|^2\,.
\end{eqnarray*}
Taking $L^1([0,1]^2)$ norm in variables $x$ and $y$ of both sides and using \eqref{eas}, \eqref{eas0}, \eqref{eas1} and  the Cauchy-Schwartz inequality gives
\[
\frac 14\int_0^1\int_0^1 \|(U_+(x,y)-U_-(x,y))e_1\|^2dxdy=\int_0^1\int_0^1 \left\|\left(\int_y^x Vds\right)e_1\right\|^2dxdy+r,\,\, |r|\lesssim \delta^{2.5}\exp(C(d+\delta))\,.
\]
Recalling the definition \eqref{okl1}, we get
\[
\left\|\left(\int_y^x Vds\right)e_1\right\|^2=\sum_{j=1}^2 (g_j(x)-g_j(y))^2
\]
so 
\[
\frac 12\int_0^1\int_0^1 \|(U_+(x,y)-U_-(x,y))e_1\|^2dxdy=4\sum_{j=1}^2\int_0^1 |g_j-\langle g_j\rangle|^2dx+r,\,\, |r|\lesssim \delta^{2.5}\exp(C(d+\delta))\,.
\]
Lemma \ref{P2le1} is proved. 
\qed

\medskip

\noindent{\bf Remark.} All statements in this subsection can be easily adjusted to any interval but the constants in the inequalities will depend on the size of that interval.
\bigskip

\subsection{Rough bound when $\widetilde \K_Q$ is small}\label{P2sec3}
\begin{Lem}
Suppose an absolutely continuous function $f$ is defined on $[0,1]$ and satisfies
\begin{eqnarray}f\in L^2[0,1],\quad  f'=l_1+l_2,\quad l_1\in L^1[0,1], \quad l_2\in L^2[0,1]\,.
\end{eqnarray}
Then, 
$
\|f\|_{L^\infty[0,1]}\le
\sqrt{\delta^2+2(\delta\tau+\epsilon(\tau+\epsilon+\delta))},
$
where $\delta=\|f\|_{L^2[0,1]},\, \epsilon=\|l_1\|_{L^1[0,1]},\, \tau=\|l_2\|_{L^2[0,1]}\,.
$
\end{Lem}
\beginpf
There is $\xi\in [0,1]$ such that $|f(\xi)|\le \delta$ and 
\[
|f(x)-f(\xi)|\le \left|\int_\xi^x f'ds\right|\le \tau+\epsilon\,.
\]
Thus, $\|f\|_{L^\infty[0,1]}\le \tau+\epsilon+\delta$. Then, writing
\[
f^2(x)-f^2(y)=2\int_y^xff'ds\,,
\]
integrating in $y$ and maximizing in $x$, we get
\[
\|f\|^2_{L^\infty[0,1]}\le \delta^2+2(\delta\tau+\epsilon(\tau+\epsilon+\delta))\,.
\]
 \qed

\medskip

Suppose $Q$ is  real-valued, symmetric matrix-function on $\R$ with zero trace and $\|Q\|_{L^2(\R)}<\infty$. Define $\Hh_Q=N^*N$, where $N: N'=JQN, N(0)=I$. Notice that $\det \int_n^{n+2} S^*\Hh_QSdx=\det \int_n^{n+2} \Hh_Qdx$ 
for every constant matrix $S\in {\rm SL}(2,\R)$. Therefore, we can apply Lemma \ref{P2le0} to each interval $[n,n+2]$ by choosing $S=N^{-1}(n)$ and get an estimate which explains how $\|Q\|_{L^2(\R)}$ controls $\widetilde \K_Q$:
\[
\widetilde \K_Q=\sum_{n\in \Z} \left(\det \int_{n}^{n+2}\Hh_Qdx-4\right)\lesssim \sum_{n\in \Z}\|Q\|^2_{L^2[n,n+2]}\exp(C\|Q\|_2)\lesssim \|Q\|^2_{L^2(\R)}\exp(C\|Q\|_2)\,.
\]
The next lemma shows that $\widetilde \K_Q$ controls the convolution of $Q$ with the exponential.

\begin{Lem}\label{P2l11}
Suppose $Q$ is  real-valued, symmetric $2 \times 2$ matrix-function on $\R$ with zero trace and entries in $L^2(\R)$. Define $\Hh_Q=N^*N$ where $N: N'=JQN, N(0)=I$ and assume that $\widetilde \K_Q<\infty$. If $O:=e^x\int_x^\infty e^{-s}Qds$, then $\|O\|_{L^\infty(\R)}\lesssim \exp(C(\|Q\|_{L^2(\R)}+\widetilde \K_Q)) \widetilde \K_Q^{{\frac 14}}
$
where  $C$ is a positive absolute constant.
\end{Lem}
\beginpf Let $R=\|Q\|_{L^2(\R)}$ and $E=\widetilde \K_Q$. We split the proof into several steps.

\noindent {\bf 1.  {Bound for a single interval $[0,1]$.}} The definitions \eqref{faw} and \eqref{faw1} imply that $\widetilde\K_{Q}^{+}\le E$.   From Theorem~1.2 and Theorem~3.2 in \cite{BD2019}, we know that $\Hh_Q$ admits the following factorization on $\R_+$: $\Hh_Q=G^*WG$ where $G$ and $W$ satisfy conditions:
\begin{eqnarray}
G'=J(v_1+v_2)G, \quad\|v_1\|_{L^1(\R_+)}\lesssim E, \quad\|v_2\|_{L^2(\R_+)}\lesssim E^{{\frac 12}}\,,\\
\det G=1,\quad v_1+v_2=(v_1+v_2)^*\,,
\end{eqnarray}
and 
\begin{eqnarray*}
W\ge 0, \quad \det W=1,\quad \|\text{tr} \,W-2\|_{L^1(\R_+)}\lesssim E\,.
\end{eqnarray*}
Since $\|\text{tr} \,W-2\|_{L^1[0,1]}\lesssim E$, we have $\|\lambda+\lambda^{-1}-2\|_{L^1[0,1]}\lesssim E$, where $\lambda$ is the largest eigenvalue of $W$. If one denotes $p=\text{tr} \,W-2=\lambda+\lambda^{-1}-2$, then 
\begin{equation}\label{P2io1}
\lambda=\frac{2+p+\sqrt{4p+p^2}}{2}\,, \quad \lambda^{-1}=\frac{2+p-\sqrt{4p+p^2}}{2}\,.
\end{equation}
In particular, that yields
\begin{equation}\label{P2p55}
\int_0^1 \|W\|\,dx\lesssim 1+E\,.
\end{equation}
The given conditions on $Q$ and \eqref{lok} yield
\[
\|N(x)\|, \|N^{-1}(x)\|\lesssim \exp(CR), \quad x\in [0,1]\,,
\]
where the second estimate follows from the first since $\det N=1$.
The Hamiltonian  $\Hh_Q=N^*N$ is absolutely continuous on $\R_+$ and 
\begin{equation}\label{tot}
0<\exp(-CR) {I}\lesssim  \Hh_Q(x)\lesssim \exp(CR) {I}
\end{equation}
 on $[0,1]$. 
We claim that $\|G(0)\|\lesssim \exp(C(R+E))$ and that $\|G^{-1}(0)\|\lesssim \exp(C(R+E))$. Indeed, if $X$ satisfies $X'=J(v_1+v_2)X$ and $X(0)=I$, then $G=XG(0)$. Moreover, given conditions on $v_1$ and $v_2$ and $\det X=1$, we have  
\begin{equation}\label{tot1}
\|X(x)\|\lesssim \exp(CE), \qquad \|X^{-1}(x)\|\lesssim \exp(CE)
\end{equation}
 uniformly on $[0,1]$. Identity $\Hh_Q=G^*(0)X^*WXG(0)$ yields
\[
(G^*(0))^{-1}\Hh_Q(G(0))^{-1}=X^*WX\,.
\]
Taking an arbitrary $\xi\in \C^2$ with $\|\xi\|_{\C^2}=1$, we get
\begin{eqnarray*}
\|G^{-1}(0)\xi\|^2\stackrel{\eqref{tot}}{\lesssim} \exp(CR) \int_0^1 \langle \Hh_Q G^{-1}(0)\xi, G^{-1}(0)\xi\rangle  dx
\hspace{3cm}
\\=  \exp(CR)\int_0^1 \langle WX\xi,X\xi\rangle dx\stackrel{\eqref{P2p55}+\eqref{tot1}}{\lesssim} \exp(C(R+E))\,,
\end{eqnarray*}
which implies $\|G^{-1}(0)\|\lesssim \exp(C(R+E))$.  We also have  $\|G(0)\|\lesssim \exp(C(R+E))$ since $\det G=1$ and the claim is proved. Finally, we have
\[
\|G(x)\|\lesssim \exp(C(R+E)), \quad \|G^{-1}(x)\|\lesssim \exp(C(R+E))
\]
for $x\in [0,1]$ since $G=XG(0)$.

Next, let us study $W$ and $W'$. Since
$W=(G^*)^{-1}N^*NG^{-1}$, one has $\|W\|\lesssim \exp(C(R+E))$ on $x\in [0,1]$. Recall that $W\ge 0$ and $\det W=1$, so
\[
\exp(-C(R+E)) {I}\lesssim W\lesssim \exp(C(R+E)) {I},\quad  x\in [0,1]\,.
\]
Since $\lambda$ is the largest eigenvalue of $W$ and $\lambda\lesssim  \exp(C(R+E))$, then \eqref{P2io1} yields $\|\lambda-1\|_{L^2[0,1]}\lesssim   E^{{\frac 12}} \exp(C(R+E))$ and $\|\lambda^{-1}-1\|_{L^2[0,1]}\lesssim  E^{{\frac 12}}\exp(C(R+E))$. Introduce $
\Upsilon=W-I\,.
$ The matrix $\Upsilon$ is unitarily equivalent to $\left( \begin{smallmatrix} \lambda-1&0\\0&1/\lambda-1\end{smallmatrix}\right)$
and that gives 
\begin{equation}\label{tot6}
\|\Upsilon\|_{L^2[0,1]}\lesssim  E^{{\frac 12}}\exp(C(R+E))\,.
\end{equation}
We need to study $\Upsilon'$, which is equal to $W'$. To do so,  notice that
\begin{eqnarray}\label{P2st1}
2N^*JQN=\Hh_Q'=G^*J(v_1+v_2)WG+G^*WJ(v_1+v_2)G+G^*W'G\,.
\end{eqnarray}
Hence,
\[
\Upsilon'=W'=F_1+F_2\,,
\]
where
 {
\[
F_1=
-J(v_1+v_2)W-WJ(v_1+v_2),\quad F_2=2(G^*)^{-1}N^*JQNG^{-1}.
\]
}
The previously obtained estimates give us
\begin{equation}\label{tot7}
 \|F_1\|_{L^1[0,1]}\lesssim   E^{{\frac 12}}\exp(C(R+E)), \quad \|F_2\|_{L^2[0,1]}\lesssim \exp(C(R+E))\,.
\end{equation}
Now, we use \eqref{tot6}, \eqref{tot7} to apply the previous lemma to each component of $\Upsilon$ to obtain 
\begin{equation}\label{tot9}
\|\Upsilon\|_{L^\infty[0,1]}\lesssim  E^{{\frac 14}}\exp(C(R+E))\,.
\end{equation} The formula \eqref{P2st1} also gives an expression for $Q$:
\[
Q=-J(H_1+H_2)\,,
\]
where
\[
H_1=0.5(N^*)^{-1}(G^*J(v_1+v_2)WG+G^*WJ(v_1+v_2)G)N^{-1}
\]
and
\[
H_2=0.5(N^*)^{-1}(G^*\Upsilon'G)N^{-1}\,.
\]
Since $\|H_1\|_{L^1[0,1]}\lesssim   E^{{\frac 12}}\exp(C(R+E))$, we have
\[
\left\|e^x\int_x^1 e^{-s}H_1ds\right\|_{L^\infty[0,1]}\lesssim  E^{{\frac 12}}\exp(C(R+E))\,.
\]
 {For smooth matrix-functions $u_1$, $u_2$, $u_3$, we have
$$
\int^{1}_{x}u_1 u'_2 u_3 \,ds = u_1u_2u_3\Bigr\rvert^{1}_{x} - \int^{1}_{x}u'_1 u_2 u_3 \,ds - \int^{1}_{x}u_1 u_2 u_3' \,ds.
$$
}
Then, 
\begin{eqnarray*}
2e^x\int_x^1 e^{-s}H_2ds=\hspace{5cm}\\e^x\left(e^{-1}(N^*(1))^{-1}G^*(1)\Upsilon(1)G(1)(N(1))^{-1}-e^{-x}(N^*(x))^{-1}G^*(x)\Upsilon(x)G(x)(N(x))^{-1}\right)\\
-e^x\int_x^1(e^{-s}(N^*(s))^{-1}G^*)'\Upsilon GN^{-1}ds-e^x \int_x^1e^{-s}(N^*(s))^{-1}G^*\Upsilon (GN^{-1})'ds\,.
\end{eqnarray*}
Since $\|(N^{-1})'\|_{L^2[0,1]}\lesssim \exp(C(R+E))$ and $\|G'\|_{L^1[0,1]}\lesssim \exp(C(R+E))$, we have
\[
\left\|e^x\int_x^1e^{-s}H_2ds\right\|_{L^\infty[0,1]}\lesssim   \|\Upsilon\|_{L^\infty[0,1]}\exp(C(R+E))\stackrel{\eqref{tot9}}{\lesssim}  E^{{\frac 14}}\exp(C(R+E))\,.
\]
Summing up, we get
\begin{equation}\label{gas1}
\left\|e^x\int_x^1e^{-s}Qds\right\|_{L^\infty[0,1]}\lesssim  E^{{\frac 14}}\exp(C(R+E))\,.
\end{equation}

\smallskip

\noindent {\bf 2. Handling all intervals $[n,n+1], n\in \Z$.} Take any $n\in \Z$. Our immediate goal is to show the bound
\begin{equation}\label{gas2}
\left\|e^x\int_x^{n+1}e^{-s}Qds\right\|_{L^\infty[n,n+1]}\lesssim  E^{{\frac 14}}\exp(C(R+E))
\end{equation}
analogous to \eqref{gas1} but written for  interval $[n,n+1]$. To this end, take the Hamiltonian $\Hh^{(n)}(x):=\Hh_Q(x+n)$ defined on $\R_+$. For the corresponding $\widetilde \K_{(n)}^+$, we get
$\widetilde \K_{(n)}^+\le E$
as follows from its definition. Since the $\widetilde \K$-characteristics of the Hamiltonians $\Hh$ and $S^*\Hh S$ are equal for every constant matrix $S\in {\rm SL}(2,\R)$, we can instead consider $\widehat \Hh^{(n)}=\widehat N^*\widehat N$ where $\widehat N'=JQ(x+n)\widehat N, \widehat N(0)=I$.  Using the arguments in step 1 for $\widehat \Hh^{(n)}$, we get \eqref{gas2}.

\smallskip

\noindent {\bf 3. Summing up.} Denote $O_n(x)=e^x\int_x^\infty e^{-s}Q\cdot \chi_{n<s<n+1}ds$ and notice that $O=\sum_{n\in \Z}O_n$. Then, since $O_n(x)=0$ for $x>n+1$ and $\|O_n(x)\|\lesssim e^{x-n}\|O_n\|_{L^\infty[n,n+1]}$ for $x<n$, we get
\[
\|O(x)\|\le \sum_{n\in \Z}\|O_n(x)\|\lesssim E^{{\frac 14}}\exp(C(R+E))\sum_{n\ge 0}e^{-n}\sim E^{{\frac 14}}\exp(C(R+E))
\]
as follows from \eqref{gas2}. That finishes the proof  {of Lemma \ref{P2l11}}. 
\qed

\subsection{ {Proof of Theorem \ref{P2c1b}}}\label{P2sec4}
Denote $E=\widetilde \K_Q$, $O=e^x\int_x^\infty e^{-s}Qds$, and recall that $\|O\|_{L^2(\R)}\sim \|Q\|_{H^{-1}(\R)}\le \|Q\|_{L^2(\R)}$.

\bigskip 

\noindent {\bf 1. Lower bound.} 
Define $\delta_n=\|O\|_{L^2[n,n+1]}$. By Lemma \ref{P2l11}, we know that $\sup_n \delta_n\lesssim  E^{{\frac 14}}\exp(C(R+E))$. Next, we apply Lemma \ref{P2le1}  to each interval $[n, n+2]$. The remainder $r_n$ in that lemma allows the estimate 
\[
r_n\lesssim (\delta_n+\delta_{n+1})^{2.5}\exp(C(\delta_n+\delta_{n+1}+R)), \quad n\in \Z\,.
\]
For each $R>0$ and  {$\eta > 0$,} we can find a positive $E_0(R,\eta)$ such that $E\in (0,E_0(R,\eta))$ implies 
that the remainder  $r_n$ is smaller than $ {\eta}(\delta_{n}^2+\delta_{n+1}^2)$ uniformly in all $n$. For example, one can take
\begin{equation}\label{raz2}
E_0(R,\eta)\sim e^{-C_\eta R},
\end{equation}
where $C_\eta $ is a sufficiently large positive number that depends on $\eta$.
Therefore, for such $E$ and  {some positive constant  $c$ independent of $\eta$}, we have
\begin{align*}
\sum_{n\in \Z} (c-\eta)\delta_n^2 \lesssim \sum_{n\in \Z} \left(\det \int_{n}^{n+2}\Hh_Q dx-4\right)\lesssim \sum_{n\in \Z} (c+\eta)\delta_n^2,
\end{align*}
where the Proposition \ref{P2p1lol} has been applied to the terms 
 {
$
\int_n^{n+2} |g_j-\langle g_j\rangle|^2\,dx
$}
in the right-hand side of \eqref{resid},  { adjusted to the interval $[n,n+2]$, to show that they are comparable to $\delta_n^2+\delta_{n+1}^2$.
Taking $\eta = c/2$, we see that
$$
E=\sum_{n\in \Z} \left(\det \int_{n}^{n+2}\Hh_Q dx-4\right)\sim \sum_{n\in \Z} \delta_n^2\sim \|O\|_{L^2(\R)}^2,
$$
in the case $E\le E_0(R,\frac{c}{2})$.} 
If $E>E_0(R,\frac{c}{2})$, one uses inequality $\|O\|_{L^2(\R)}\lesssim R$ to get
\begin{equation}\label{raz3}
e^{-CR}\|O\|^2_{L^2(\R)}\lesssim \frac{E_0(R,\frac{c}{2})}{1+R^2}\|O\|^2_{L^2(\R)}\lesssim E\,,
\end{equation}
which holds for some positive absolute constant $C$ due to \eqref{raz2}. That provides the required lower bound.

\bigskip

\noindent{\bf 2. Upper bound.} Let
$
\delta_n=\|O\|_{L^2[n,n+1]}\,.
$ For given value of $R$,   apply Lemma~\ref{P2le1} and Proposition~\ref{P2p1lol} to each interval $[n,n+2]$. That gives
\[
E\lesssim \sum_{n\in \Z} \delta^2_ne^{C(R+\delta_n)}
\]
with an absolute constant $C$. Since $\sum_{n\in \Z}\delta_n^2\sim \|q\|_{H^{-1}(\R)}^2$ and $\|q\|_{H^{-1}(\R)}\lesssim R$, one has
\[
E\lesssim \|q\|_{H^{-1}(\R)}^2e^{C(R+\|q\|_{H^{-1}(\R)})}\lesssim \|q\|_{H^{-1}(\R)}^2e^{C_2R}\,.
\]
\qed
\bigskip

\bigskip

\section{Appendix}
 {Here we collect} some auxiliary results used in the main text.  

\medskip

{\bf 1.} We start with an example that shows that the scattering transform is not injective when defined on $q\in L^2(\R)$. This is an analog of Lemma 17 in \cite{TT}. 

\medskip

\begin{Ex}\label{ex1}
There exist potentials $q_1, q_2 \in L^2(\R)$ such that $q_1 \neq q_2$ in $L^2(\R)$ but we have $\rc_{q_1} = \rc_{q_2}$ a.e.~on $\R$ for their reflection coefficients. In other words, the scattering transform $q \mapsto \rc_{q}$ is not injective on $L^2(\R)$.
\end{Ex}
\beginpf Let us consider 
$$
\fa_1^{+} = 1, \quad \fb_1^{+} = 0,  \quad \fa_1^{-} = \fa, \quad \fb_1^{-} = \fb,
$$
and 
$$
\fa_2^{+} = \fa, \quad \fb_2^{+} = \fb, \quad
\fa_2^{-} = 1, \quad \fb_2^{-} = 0,  
$$
where $\fa = 1+i/x$, $\fb = i/x$. Note that 
$$
\int_{\R}\log(1 - |s_{k}^{\pm}(x)|^2)\,dx > -\infty, \qquad s_{k}^{\pm} = \frac{\fb^{\pm}}{\fa_{k}^{\pm}}, \qquad k=1,2.
$$
Theorem 12.11 in \cite{Den06} says that for every contractive analytic function $s$ on $\C_+$ whose boundary values on $\R$ satisfy $\log (1-|s|^2) \in L^1(\R)$ there exists a unique coefficient $A \in  L^2(\R_+)$ such that $s = \lim_{\xi\to+\infty}\frac{\fB(\xi,\lambda)}{\fA(\xi,\lambda)}$, $\lambda \in \C_+$ for the continuous Wall polynomials generated by $A$. Moreover, we have 
\begin{equation}\label{eq11}
2\pi\|A\|_{L^2(\R_+)}^{2} = \|\log(1-|s|^2)\|_{L^1(\R_+)}.
\end{equation}
Applying this result, we see that there exist functions $A_{1}^{\pm}, A_{2}^{\pm} \in L^2(\R_+)$ such that $\fa^{\pm}_{1,2}$, $\fb_{1,2}^{\pm}$ are the limits of their continuous Wall polynomials. Now define potentials $q_{1, 2} \in L^2(\R)$ by relations
$$
A_{1,2}^+(\xi) = -\ov{q_{1,2}(\xi/2)}/2, \qquad 
A_{1,2}^{-}(\xi) = q_{1,2}(-\xi/2)/2, \qquad \xi \in \R_+.
$$
From Proposition \ref{p21}, we conclude that the coefficients $a_{1,2}$, $b_{1,2}$ for these potentials satisfy
$$a_1 = \fa = a_2, \qquad b_1 = -\fb = \ov{\fb}= b_2,
$$
on $\R\setminus \{0\}$. Hence, $\rc_{q_1} = \rc_{q_2}$ on $\R\setminus \{0\}$. On the other hand, we have $A_{1}^{+} = 0$ and $A_{2}^{-} = 0$ by construction. It follows that $\supp q_1 \subset (-\infty, 0]$ and $\supp q_2 \subset [0, +\infty)$. Since $q_1$ and $q_2$ are nonzero (they have a nonzero $L^2(\R)$-norm as follows from \eqref{eq11}), that yields $q_1 \neq q_2$ in $L^2(\R)$. \qed

\medskip

{\bf 2.} Next, we outline how to prove  that the  spectral representation for the Dirac operator $\Di_{Q}$, defined by relation \eqref{eq18}, is given by the Weyl-Titchmarsh transform \eqref{wt}. To this end, we will use the corresponding result for canonical Hamiltonian systems proved in \cite{Remlingb}.

\medskip

At first, we note that if $\Hh_{Q} = N_Q^* N_Q$ is the Hamiltonian from Theorem \ref{t15}, then $\det \Hh_{Q} = 1$ on $\R$ and the operator $V: X \mapsto N_{Q}^{-1} X$ is unitary from $L^2(\R, \C^2)$ onto the Hilbert space 
$$
L^2(\Hh_Q) = \Bigl\{Y: \R \to \C^2: \; \|Y\|_{L^2(\Hh_Q, \R)}^{2} = \int_{\R}\langle\Hh_Q(\xi)Y(\xi), Y(\xi)\rangle_{\C^2}\,d\xi < \infty\Bigr\}.
$$ 
Moreover, $V \Di_{Q} V^{-1}$ coincides with the operator $\Di_{\Hh_{Q}}: Y \mapsto \Hh^{-1}J Y'$ of the canonical Hamiltonian system generated by the Hamiltonian $\Hh_{Q}$. Thus, the operator $\Di_{Q}$ on $L^2(\R, \C^2)$ is unitary equivalent to the operator $\Di_{\Hh_{Q}}$ on $L^2(\Hh_Q)$. Let $\widetilde M$ be the solution of Cauchy problem 
\begin{equation}\label{eq21}
J\widetilde M'(\xi, z) = z \Hh_{Q}(\xi)\widetilde M(\xi, z), \qquad 
\widetilde M(0, z) = \idm,
\end{equation}
where $z \in \C$, $\xi \in \R$, and the differentiation is taken with respect to $\xi \in \R$. The Weyl-Titchmarsh transform for $\Di_{\Hh_Q}$ is defined by
$$
\F_{\Di_{\Hh_{Q}}}: Y \mapsto \frac{1}{\sqrt{\pi}}\int_{\R}\widetilde M(\xi,\lambda)^* \Hh_{Q}(\xi) Y(\xi)\,d\xi
$$
on a dense subset of $L^2(\Hh_Q)$ of smooth compactly supported functions. This operator is unitary from $L^2(\Hh_Q)$ onto the space $L^2(\rho)$ defined in the same way as at the beginning of Section \ref{s3}. Specifically, we let $m_{\pm}$ be the half-line Weyl functions of $\Hh_Q$ and define $\rho$ as the representing measure for the matrix-valued Herglotz function $m$ in \eqref{eq22}. It was proved in Theorem 3.21 in \cite{Remlingb} that $\F_{\Di_{\Hh_{Q}}}\Di_{\Hh_{Q}}\F_{\Di_{\Hh_{Q}}}^{-1}$ coincides with the operator of multiplication by the independent variable in $L^2(\rho)$. We also have 
$$
\F_{\Di_{\Hh_{Q}}}(VX) = \F_{\Di_{{Q}}} X, \qquad X \in L^2(\R, \C^2).
$$
Thus, we only need to check that the Weyl functions $m_{\pm}$ used in Section \ref{s3} coincide with the half-line Weyl functions of the Hamiltonian $\Hh_{Q}$. For the $\R_+$-Weyl functions this follows from Lemma~\ref{l6} below. Comparing the formulas for $A^+$, $A^-$ in the beginning of Section \ref{s3}, we see that the Weyl function $m_-$ for $\Di_{Q}$ corresponds to the Weyl function $m_+$ for $\Di_{\widetilde Q}$ where $\widetilde Q(\xi) = \sigma_{3} Q(-\xi) \sigma_{3}$. Similarly, in the setting of canonical Hamiltonian systems, the Weyl function $m_-$ for $\Di_{\Hh_{Q}}$ coincides with the Weyl function $m_+$ of $\Di_{\widetilde\Hh_{Q}}$, $\widetilde \Hh_{Q}(\xi) = \sigma_{3}\Hh(-\xi)\sigma_{3}$. Therefore, the statement for $A^-$  follows from Lemma \ref{l6} below and from the relation $\widetilde \Hh_{Q} = \sigma_3 \Hh_{Q} \sigma_3 = (\sigma_3 N_Q^* \sigma_3) (\sigma_3 N_Q \sigma_3) = \Hh_{\widetilde Q}$.    

\begin{Lem}\label{l6}
Let $q \in L^2(\R_+)$. Define
$$
Q(\xi) = \begin{pmatrix}-\Im q(\xi) & -\Re q(\xi) \\ -\Re q(\xi) & \Im q(\xi)\end{pmatrix}, \qquad A(\xi) = -\ov{q(\xi/2)}/2, \qquad \xi \in \R_+.
$$
Let $N_Q$ be defined by $JN'_Q(\xi,\lambda) + Q(\xi) N_Q(\xi,\lambda) = \lambda N_Q(\xi,\lambda)$, $N_Q(0,\lambda) = \idm$. Consider the Hamiltonian $\Hh_{Q} = N_Q^*(\xi,0)N_Q(\xi,0)$ on $\R_+$ and let $\widetilde M = \left(\begin{smallmatrix}
\widetilde M_{11} & \widetilde M_{12}\\
\widetilde M_{21} & \widetilde M_{22}
\end{smallmatrix}\right)$ be defined by $J\widetilde M'(\xi, z) = z \Hh_{Q}(\xi)\widetilde M(\xi, z)$, 
$\widetilde M(0, z) = \idm$. Let, finally, $P$, $P_*$, $\widehat P$, $\widehat P_{*}$ be the solutions to Krein systems \eqref{eq7}, \eqref{eq8} for the coefficient $A$ on $\R_+$. Then,
\begin{equation}\label{eq24}
\lim_{\xi \to +\infty} \frac{\widetilde M_{22}(\xi,z)}{\widetilde M_{21}(\xi,z)} =
\lim_{\xi \to +\infty} \frac{(N_Q)_{22}(\xi,z)}{(N_Q)_{21}(\xi,z)} 
= \lim_{\xi \to +\infty} i\frac{\widehat P_{*}(\xi,z)}{P_{*}(\xi,z)}, \qquad z \in \C_+.
\end{equation}
In other words, the function $m_{+}$ in \eqref{eq23} is the half-line Weyl function for the operators $\Di_{\Hh_Q}$, $\Di_{Q}$.
\end{Lem}
\beginpf The formula
$$
\lim_{\xi \to +\infty} \frac{\widetilde M_{22}(\xi,z)}{\widetilde M_{21}(\xi,z)} =
\lim_{\xi \to +\infty} \frac{(N_Q)_{22}(\xi,z)}{(N_Q)_{21}(\xi,z)}
$$
for $\Di_{Q}$ and $\Di_{\Hh_{Q}}$ is well-known and can be derived from the analysis of Weyl circles  by using identity  $N_Q(\xi,\lambda) = N_Q(\xi,0) \widetilde M(\xi, \lambda)$ and the invariance of Weyl circles under transforms generated by $J$-unitary matrices (in our setting, the $J$-unitary matrix is $N_Q(\xi, 0)$: we have $N_Q^*(\xi, 0) J N_Q(\xi, 0) = J$ on $\R$). See, e.g., \cite{BLY1} or Section~8 in \cite{Romanov} for more details on Weyl circles for canonical Hamiltonian systems. Thus, we focus on the second identity in \eqref{eq24} and define
$$
X(\xi,z)=
e^{-i\xi z}
\begin{pmatrix}
\displaystyle \frac{P(2\xi,z)+P_*(2\xi,z)}{2} & \displaystyle\frac{\widehat P(2\xi,z) - \widehat P_*(2\xi,z)}{2i}\\
\displaystyle \frac{P_*(2\xi,z)-P(2\xi,z)}{2i} & \displaystyle\frac{\widehat P(2\xi,z) + \widehat P_*(2\xi,z)}{2}
\end{pmatrix}, \qquad \xi \in \R, \quad z \in \C\,.
$$
Differentiating, one obtains $JX' + QX = z X$, $X(0, z) = \idm$. It follows that $X(\xi,z) = N_Q(\xi, z)$. In particular, we have 
$$(N_Q)_{22} = e^{-i\xi z}\frac{\widehat P(2\xi,z) + \widehat P_*(2\xi,z)}{2}, \qquad (N_Q)_{21} = e^{-i\xi z}\frac{P_*(2\xi,z)-P(2\xi,z)}{2i}.$$ 
Since $P(\xi, z) \to 0$, $P_*(\xi, z) \to \Pi(z) \neq 0$ as $\xi \to +\infty$ (see Theorem 12.1 in \cite{Den06}), and analogous relations hold for $\widehat P$ and $\widehat P_*$, we have
$$
\lim_{\xi \to +\infty} \frac{(N_Q)_{22}(\xi,z)}{(N_Q)_{21}(\xi,z)} 
= 
\lim_{\xi \to +\infty} i\frac{\widehat P_{*}(\xi,z)}{P_{*}(\xi,z)}, \qquad z \in \C_+.
$$  
The lemma is proved. \qed

\bigskip

{\bf 3.} Lemma \ref{l6} and some known results for canonical systems can be used to show that weak convergence of potentials  {of} the Dirac operator implies  convergence of the corresponding Weyl functions.
\begin{Lem}\label{l7}
Suppose $\{q_{\ell}\}_{\ell > 0}$ is a bounded sequence  in $L^2(\R_+)$ which converges to zero weakly. Let $Q_{\ell}$ be the associated matrix-functions  defined as in Lemma \ref{l6}. Then, the sequence of corresponding Weyl functions
$\{m_{\ell, +}\}$
converges to $i$ locally uniformly in $\C_+$ when $\ell\to +\infty$.
\end{Lem}
\beginpf For $\ell>0$, denote by $\Hh_{Q_\ell}$ the Hamiltonian generated by $Q_\ell$ as in Lemma \ref{l6}. Then, $m_{\ell, +}$ is the Weyl function for the half-line operators $\Di_{\Hh_{Q_\ell}}$ and $\Di_{Q_{\ell}}$. Since $\sup_{\ell > 0}\|q_\ell\|_{L^2(\R_+)} < \infty$ and  $q_\ell$ converge to zero weakly in $L^2(\R_+)$ as $\ell \to +\infty$, the Hamiltonians $\Hh_{Q_\ell}$ tend to the identity matrix $\Hh_{0} = \idm$ uniformly on compact subsets on $\R_+$. Then, their Weyl functions $m_{+,\ell}$ tend to the Weyl function $m_+ = i$ of the Hamiltonian $\Hh_{0}$ locally uniformly in $\C_+$ by Theorem $5.7$ $(b)$  in \cite{Remlingb}. \qed

\bibliographystyle{plain} 
\bibliography{bibfile}

\end{document}